\documentclass[preprint]{elsarticle}

\usepackage{adjustbox}
\usepackage{algorithmic}
\usepackage{amsmath,amssymb,amsfonts}
\usepackage{appendix}
\usepackage{array}
\usepackage{bm}
\usepackage{booktabs}
\usepackage{diagbox}
\usepackage{dsfont}
\usepackage{graphicx}
\usepackage{hyperref}
\usepackage{mathtools}
\usepackage{multirow}
\usepackage{siunitx}
\usepackage{soul}
\usepackage{stmaryrd}
\usepackage[caption=false]{subfig}
\usepackage{textcomp}
\usepackage{xcolor}

\SetSymbolFont{stmry}{bold}{U}{stmry}{m}{n}

\newcommand{\norm}[1]{\left\lVert#1\right\rVert}
\newcommand{\abs}[1]{\left\lvert#1\right\rvert}
\newcommand{\nint}[1]{\llbracket#1\rrbracket}
\newcommand{\diff}{\mathop{}\!\mathrm{d}}

\def\*#1{\mathbf{#1}}

\newcommand{\bxini}{\mathbf{\overline{x}}}
\newcommand{\mD}{\mathcal{D}}
\newcommand{\mJ}{\mathcal{J}}
\newcommand{\mK}{\mathcal{K}}
\newcommand{\mL}{\mathcal{L}}
\newcommand{\mM}{\mathcal{M}}
\newcommand{\mR}{\mathcal{R}}
\newcommand{\mX}{\mathcal{X}}
\newcommand{\NN}{\mathbb{N}}
\newcommand{\RR}{\mathbb{R}}
\newcommand{\Sym}{\mathbb{S}}

\DeclareMathOperator{\argmin}{argmin}
\DeclareMathOperator{\degree}{degree}
\DeclareMathOperator{\Card}{Card}
\DeclareMathOperator{\Diag}{Diag}
\DeclareMathOperator{\Id}{\mathrm{Id}}

\begin{document}


\begin{frontmatter}

  \title{Sparse Signal Reconstruction for Nonlinear Models via Piecewise
    Rational Optimization}
  
  \author[CVN]{Arthur Marmin\corref{my_corres_author}}
  \cortext[my_corres_author]{Corresponding author}
  \ead{arthur.marmin@centralesupelec.fr}
  
  \author[TSP]{Marc Castella}
  \ead{marc.castella@telecom-sudparis.eu}
  
  \author[CVN]{Jean-Christophe Pesquet}
  \ead{jean-christophe.pesquet@centralesupelec.fr}
  
  \author[IFPEN]{Laurent Duval}
  \ead{laurent.duval@ifpen.fr}

  \address[CVN]{University Paris-Saclay, CentraleSup{\'e}lec Center for Visual
    Computing, Inria, \\
    9 Rue Joliot Curie, 91190 Gif-sur-Yvette, France}
  \address[TSP]{SAMOVAR, CNRS, T{\'e}l{\'e}com SudParis, Institut Polytechnique
    de Paris, \\
    91011 Evry Cedex, France}
  \address[IFPEN]{ESIEE Paris, University Paris-Est, LIGM, Noisy-le-Grand, \\
    and IFP Energies nouvelles, Rueil-Malmaison, France.}
  
  
  \begin{abstract}
    We propose a method to reconstruct sparse signals degraded by a nonlinear
    distortion and acquired at a limited sampling rate.
    Our method formulates the reconstruction problem as a nonconvex minimization
    of the sum of a data fitting term and a penalization term.
    In contrast with most previous works which settle for approximated local
    solutions, we seek for a global solution to the obtained challenging
    nonconvex problem.
    Our global approach relies on the so-called Lasserre relaxation of
    polynomial optimization.

    We here specifically include in our approach the case of piecewise rational
    functions, which makes it possible to address a wide class of nonconvex
    exact and continuous relaxations of the $\ell_{0}$ penalization function.
    Additionally, we study the complexity of the optimization problem.
    It is shown how to use the structure of the problem to lighten the
    computational burden efficiently.
    Finally, numerical simulations illustrate the benefits of our method in
    terms of both global optimality and signal reconstruction.
  \end{abstract}
  
  

  \nonumnote{Preliminary versions of this work were presented
    in~\cite{Marmin_A_2018_p-eusipco_signal_rsndo}
    and in~\cite{Marmin_A_2019_p-icassp_how_gsnopialp}}

\end{frontmatter}


\section{Introduction}

Sparse signals, i.e.\ signals composed of a few spikes, are of particular
interest.
They either occur naturally in many areas or emerge after
sparsifying transformations such as time-frequency or wavelet
decompositions~\cite{Gauthier_J_2009_j-ieee-trans-sig-proc_optimization_socfb,
  Pham_M_2014_j-ieee-trans-sig-proc_primal-dual_pastbafasmr}.
However, accurate data acquisition of sparse signals from real-world
measurements remains an open challenge.
The difficulty of the problem is further increased when acquiring data at a
reduced rate.
This is however an important practical situation, since it permits faster
acquisitions for high-throughput experiments and analysis.

A common approach to recover the original signal from the observations is
first to define a well-chosen criterion and then to minimize it.
The criterion is often composed of two terms: a fit function depending on the
investigated model as well as the observations, and a (possibly composite)
regularization term that allows good estimates to be selected among those
consistent with the data~\cite{Chaux_C_2007_j-inv-prob_variational_ffip}.
However, few methods today are able to deal with nonlinear models and to
globally optimize sparsity promoting criteria.
Indeed, integrating any of these two properties in the criterion often yields an
intricate optimization problem that is difficult to solve.

Thence, to deal with nonlinear effects, linearization techniques are often used
since the vast majority of available methods only apply to linear
models~\cite{Tibshirani_R_1996_j-r-stat-soc-b-stat-meth_regression_ssl,
  Blumensath_T_2008_j-four-anal-appl_iterative_tsa,
  Soubies_E_2015_j-siam-j-imaging-sci_continuous_elplsrp} or to models with
weaker linearity assumptions~\cite{Schetzen_M_2010_book-nonlinear_smavwp,
  Dobigeon_N_2014_j-ieee-sig-proc-mag_nonlinear_uhima,
  Deville_Y_2015_book_overview_bssmlqpnm}.
On the other hand, the standard approach to promote sparse solutions consists in
adding an $\ell_{0}$ penalization to a data-fit cost function which leads to NP
hard optimization
problems~\cite{Nikolova_M_2013_j-siam-j-imaging-sci_description_mlsrlnugm,
  Bourguignon_S_2016_j-ieee-trans-sig-proc_exact_sapmipfcp}.
Consequently, several surrogates to the $\ell_{0}$ penalization have been
suggested, the simplest one being the $\ell_{1}$ norm.
The latter has the enjoyable property of being convex, which simplifies the
optimization task~\cite{Combettes_P_2008_j-siam-j-optim_proximal_tamob,
  Combettes_P_2011_book_fixedpoint_aipse}, but it also strongly penalizes high
values of the variables and thus introduces a bias in the solutions.
Albeit providing good results, the nonconvex Geman-McClure
function~\cite{Castella_M_2015_p-ieee-camsap_optimization_gmlcssd} also tends to
introduce bias.
Therefore further relaxations of $\ell_{0}$ function have been
investigated~\cite{Fan_J_2001_j-amer-stat-asso_variable_snplop}.
A major drawback is that those relaxations are nonconvex and result in
optimization problems which are difficult to solve globally in the sense that
currently available algorithms only converge to local solutions and therefore
may be highly dependent on their
initialization~\cite{Fan_J_2001_j-amer-stat-asso_variable_snplop,
  Ochs_P_2015_j-siam-j-imaging-sci_iteratively_r_a_nnocv,
  Candes_E_2008_j-four-anal-appl_enhancing_srlm,
  Breheny_P_2011_j-ann-appl-stat_coordinate_danprabfs,
  Patrascu_A_2015_j-ieee-trans-automat-ctrl_random_cdmlrco,
  Blumensath_T_2008_j-four-anal-appl_iterative_tsa,
  Selesnick_I_2017_j-ieee-trans-sig-proc_sparse_rca}.

In the case of a linear model, a first approach for ensuring global convergence
of an exact relaxation of the $\ell_{0}$ function has been proposed
in~\cite{Bourguignon_S_2016_j-ieee-trans-sig-proc_exact_sapmipfcp} and is based on
mixed-integer programming.
This work proposes a different approach grounded on the global minimization of
the broad class of piecewise rational functions under polynomial constraints.
Based on it, we propose a novel recovery method for sparse signals from
subsampled observations obtained through a noisy model involving nonlinear
functions.
More precisely, we show that the fit function and the regularization term can be
modeled as piecewise rational functions.
Fortunately, many well-known good approximations to the $\ell_{0}$ penalization
satisfy the latter property~\cite{Zhang_T_2010_j-mach-learn-res_analysis_mcrsr,
  Fan_J_2001_j-amer-stat-asso_variable_snplop,
  Zhang_C_2010_j-ann-appl-stat_nearly_uvsumcp,
  Artina_M_2013_j-siam-j-optim_linearly_cnnm,
  Jezierska_A_2011_book_fast_stcpqcsm,
  Soubies_E_2015_j-siam-j-imaging-sci_continuous_elplsrp,
  Marmin_A_2019_p-icassp_how_gsnopialp}.
Moreover, various nonlinear degradations, such as saturation, can be modeled
with rational functions.
Hence, several criteria of interest for reconstructing sparse signals which have
been nonlinearly degraded can be modeled, or faithfully approximated, as
piecewise rational.
We then reformulate the corresponding piecewise rational optimization problem
as the minimization of a sum of rational functions, for which the recent
framework of Lasserre's hierarchy~\cite{Lasserre_J_2001_j-siam-j-optim_global_oppm}
can be applied.
This framework relaxes a polynomial optimization problem into a hierarchy of
convex semi-definite programming (SDP) problems whose solutions converge to a
global solution to the initial polynomial problem.
SDP problems are playing an important role in our methodology, however, solving
large dimensional SDP problems remains nowadays an open challenge.
Therefore, we study the overall complexity of the SDP relaxations and show how
to reduce it efficiently in several ways.
We especially emphasize the benefit of subsampling.
Our contribution is twofold:
\begin{itemize}
\item First, we investigate a wide range of continuous approximations to the
  $\ell_{0}$ penalty and we extend the framework of Lasserre's hierarchy to
  piecewise rational functions in order to minimize the resulting nonconvex
  criterion.
  Unlike standard approaches, we are able to establish theoretical guarantees on
  the global optimum of the original optimization problem.
  In particular, we provide a unified view of our previous
  works~\cite{Marmin_A_2018_p-eusipco_signal_rsndo,
    Castella_M_2019_j-ieee-trans-sig-proc_rational_onral0p,
    Marmin_A_2019_p-icassp_how_gsnopialp}.
  The framework and the nonlinear observation model have first been proposed
  in~\cite{Castella_M_2019_j-ieee-trans-sig-proc_rational_onral0p} while the
  subsampling has been introduced in~\cite{Marmin_A_2018_p-eusipco_signal_rsndo}.
  However, in both~\cite{Castella_M_2019_j-ieee-trans-sig-proc_rational_onral0p,
    Marmin_A_2018_p-eusipco_signal_rsndo}, the regularizer was restricted to a
  Geman-Mcclure potential.
  We propose here to use a much richer class of regularizers which was introduced
  in~\cite{Marmin_A_2019_p-icassp_how_gsnopialp} but only for a simple linear
  model.
\item Second, through a complexity analysis and extensive simulations, we show
  how the structure of the problem and subsampling allow us to alleviate the
  computational burden of the original Lasserre's framework.
  Our approach can be successfully applied to signal processing and compressed
  sensing problems as illustrated by the provided example inspired by the
  acquisition of signals in gas chromatography.
\end{itemize}

Our article is organized as follows: Section~\ref{sec:motivation} introduces our
model and criterion while Section~\ref{sec:l0_approx} presents the class of
approximations to the $\ell_{0}$ penalty we consider before reformulating the
minimization of our criterion as a rational optimization problem.
Section~\ref{sec:solve} details how to solve such optimization problem by
leveraging its inherent structure.
Section~\ref{sec:sdp_compl} first studies the complexity of the obtained SDP
problems before explaining how to decrease it efficiently.
Section~\ref{sec:simulations} presents numerical simulations in order to
validate our method.
Section~\ref{sec:conclusion} concludes our work.

We introduce the following notation: $\ast$ is the convolution operator, for any
nonnegative integers $n$ and $k$, $\Sym^{n}$ (resp.\ $\Sym_{+}^{n}$) is the set of
$n \times n$ real symmetric (resp.\ symmetric positive semi-definite) matrices,
$\binom{n}{k}$ is the binomial coefficient ``among $n$ choose $k$'', $\lfloor \cdot \rfloor$
(resp. $\lceil \cdot \rceil$) is the greatest (resp.\ smallest) integer lower (resp.\ greater)
than its argument, $\abs{\bm{\alpha}}=\alpha_{1}+\cdots+\alpha_{n}$ denotes the absolute value of a
multi-index $\bm{\alpha}=(\alpha_{1},\dots,\alpha_{n})$ of size $n$, and $\NN_{t}^{n}$ is the subset of
multi-indices whose absolute value is less than or equal to $t$.
The superscript $\mbox{}^\top$ indicates the transpose of a matrix.
For a given set $\mX$, $\mathds{1}_{\{\cdot \in \mX\}}$ is the characteristic
function of $\mX$ with $\mathds{1}_{\{x \in \mX\}}=1$ if $x$ is in $\mX$ and
$0$ otherwise.
For a given polynomial $p$, we define the following operator
\begin{equation}
  \label{eq:defdemideg}
  \mathrm{d}_{p} = \left\lceil\frac{\degree p}{2}\right\rceil \, ,
\end{equation}
and we denote by $\*p$ a vector composed of the coefficients corresponding
to monomials in $p$ up to the total degree of $p$.


\section{Observation and signal model}
\label{sec:motivation}

\subsection{Our observation model}

We consider the reconstruction of an unknown discrete-time sparse signal
$\bxini$ of length $T$.
The measurement process deteriorates $\bxini$ in the following way: the peaks it
contains are enlarged and the sensors introduce a saturation effect.
As common in the literature, these degradations are modeled respectively by a
convolution with a finite impulse response filter and by a memoryless nonlinear
function $\Phi$.
The filter coefficients are given by a vector $\*h$ of length $L$.
Finally, a noise is superimposed, which is modeled by an additive vector term
$\*w$ with samples drawn from an i.i.d.\@ zero-mean Gaussian distribution.

An important feature of our model is its ability to deal with subsampling of the
measured signal during the acquisition.
As in many applications such as chromatography and spectroscopy, the physical
limitations may allow only subsampled data acquisition, we introduce a
decimation operator $D$.
Interestingly, we will see that our approach is applicable in this context and
allows one to use well-suited penalization terms to promote sparsity.
Defining the observation vector $\*y$ of size $U$ after subsampling, the
corresponding modeling equation finally reads
\begin{equation}
  \label{eq:model}
  \*y = D\big(\Phi(\*h \ast \bxini) + \*w\big) \, .
\end{equation}
Model~\eqref{eq:model} can emulate narrow-peak signals from gas chromatography
experiments~\cite{Vendeuvre_C_2005_j-chromato-a_characterization_mdctdgcgcgcpapvsamd,
  Vendeuvre_C_2007_j-ogst_comprehensive_tdgcdcpp}.
In this case, the filter $\*h$ has a discretized Gaussian shape.
This choice arises from traditional stochastic or plate modeling, representing
a Galton-Hennequin bell distribution~\cite[Chapter 3]{Felinger_A_1998_book_data_aspc}.
Peak saturation is also modeled, which cannot be done in standard practice in
analytical chemistry.
According to~\cite{Kalambet_Y_2018_j-chemometr-intell-lab-syst_comparison_ircvncp},
filter lengths $L$ from $3$ to $9$ samples may suffice for a relatively accurate
estimation of the peak area, a quantity related to the concentration of a
particular molecule.

We will be interested in regular decimation patterns $D_{\delta}$ where all the
elements indexed with a multiple of an integer $\delta$ are deleted, namely
\begin{equation}
  D_{\delta}\big({(s_{t})}_{t \in \nint{1,T}}\big) = {(s_{\Delta(u,\delta)})}_{u \in \nint{1,U}}
  \, ,
\end{equation}
where
\begin{equation}
  \label{eq:Udelta}
  U = T- \lfloor T /\delta \rfloor
\end{equation}
and $\Delta$ is defined as
\[
  (\forall u \in \nint{1,U}) \quad
  \Delta(u,\delta) = u + \left\lfloor \frac{u-1}{\delta-1} \right\rfloor
  \, .
\]
We denote by $D_{\infty}$ the identity operator that preserves the entire signal.
Let us illustrate the two decimation patterns $D_{2}$ and $D_{4}$ on the example
vector $\*s={[s_{1},s_{2},s_{3},s_{4},s_{5},s_{6},s_{7},s_{8}]}^{\top}$
\begin{align*}
  \*s &\xmapsto{D_{2}} {[s_{1},s_{3},s_{5},s_{7}]}^{\top}
  = {(s_{\Delta(u,2)})}_{u \in \nint{1,4}} \\
  \*s &\xmapsto{D_{4}} {[s_{1},s_{2},s_{3},s_{5},s_{6},s_{7}]}^{\top}
  = {(s_{\Delta(u,4)})}_{u \in \nint{1,6}} \, .
\end{align*}
The smaller parameter $\delta$, the higher the decimation and harder the
reconstruction of the signal $\bxini$.

To estimate the original signal $\bxini$, we minimize a penalized criterion
$\mJ$ composed of two terms:
\begin{equation}
  \label{eq:init_crit}
  (\forall \*x \in \RR^{T}) \quad \mJ(\*x) = f_{\*y}(\*x) + \mR_{\lambda}(\*x) \, .
\end{equation}
The first one $f_{\*y}$ is a fit measure with respect to the acquired
measurements $\*y$ while the second one $\mR_{\lambda}$ is a regularization term which
will be discussed next in more detail in Section~\ref{ssec:ex_reg}.

As a fit function, we choose the standard least-squares error between $\*y$ and
the output of the noiseless model for a given estimate $\*x$ of the original
signal $\bxini$
\[
  (\forall \*x \in \RR^{T}) \quad
  f_{\*y}(\*x) = \norm{\*y-D_{\delta}(\Phi(\*h \ast \*x))}_{2}^{2}
  = \norm{\*y-D_{\delta}(\Phi(\*H\*x))}_{2}^{2}\, ,
\]
where $\*H$ is a Toeplitz band matrix corresponding to the convolution with
$\*h$.
Because of the transformation $\Phi$, the fit function $f_{\*y}$ is possibly
nonconvex.
This is in contrast with more classical linear models in which the fit function
reduces to the quadratic function $\* x \longmapsto \norm{\*y-D_{\delta}(\*H\*x)}_{2}^{2}$.
In our approach, other fit functions $f_{\*y}$ can be chosen to model different
problems as long as they are rational.
In the following, the nonlinear function $\Phi$ is assumed to be rational and to
act component-wise.
Setting the components of $\*x$ with nonpositive index to be identically zero in
order to unclutter notation, $f_{\*y}$ hence reads as a sum of rational
functions
\[
  f_{\*y}(\*x)
  = \sum_{u=1}^{U} \underbrace{{\left(y_{u}-\Phi\Big(\sum_{l=1}^{L}h_{l}x_{\Delta(u,\delta)-l+1}\Big)\right)}^{2}}
  _{\mbox{$g_{u}(x_{\Delta(u,\delta)-L+1},\ldots,x_{\Delta(u,\delta)})$}}
  \, ,
\]
where ${(g_u)}_{u \in \nint{1,U}}$ are rational functions in $L$ variables.


\subsection{Properties of the original signal and examples of
  \texorpdfstring{$\ell_{0}$}{l0} approximations}
\label{ssec:ex_reg}

The unknown original signal $\bxini$ sought by the reconstruction method is
assumed to be sparse.
In other words, it comprises only few peaks and many of its components are zero.
Following this assumption, the second term $\mR_{\lambda}$ in~\eqref{eq:init_crit} is
a sparsity-promoting penalization weighted by a positive parameter $\lambda$.
Ideally, we would like $\mR_{\lambda}$ to be the sparsity measure $\lambda\ell_{0}$ (where $\ell_0$
counts the number of nonzero elements) but, in order to derive computationally
efficient optimization techniques, a suitable separable approximation is
substituted for it, which reads
\begin{equation}
  \label{eq:sep_reg}
  (\forall \*x = {(x_{t})}_{t\in\nint{1,T}} \in \RR^{T}) \quad 
  \mR_{\lambda}(\*x) = \sum_{t=1}^{T} \Psi_{\lambda} (x_{t}) \, .
\end{equation}
Common approaches consist in using either convex functions $\Psi_{\lambda}$ such as the
$\ell_{1}$ norm, or nonconvex ones that still maintain the convexity of the overall
criterion~\cite{Selesnick_I_2017_j-ieee-trans-sig-proc_sparse_rca}.
However, a good approximation $\Psi_{\lambda}:\RR\rightarrow\RR$ to the $\ell_{0}$ function requires
the following three properties~\cite{Fan_J_2001_j-amer-stat-asso_variable_snplop}
leading to nonconvex criteria: unbiasedness for large values, sparsity to
reduce the complexity of the model by setting small values to zero, and
continuity to ensure the stability of the model.
In contrast with~\cite{Castella_M_2019_j-ieee-trans-sig-proc_rational_onral0p,
  Marmin_A_2018_p-eusipco_signal_rsndo} where the Geman-McClure nonconvex $\ell_{0}$
approximation was used and introduced bias in the estimate, we propose here a
much wider class of piecewise rational function approximations that satisfy the
three mentioned properties.
Those approximations extend significantly our previous work to settings of more
practical interest.

\begin{figure}[htbp]
  \centerline{\includegraphics[width=0.7\textwidth]{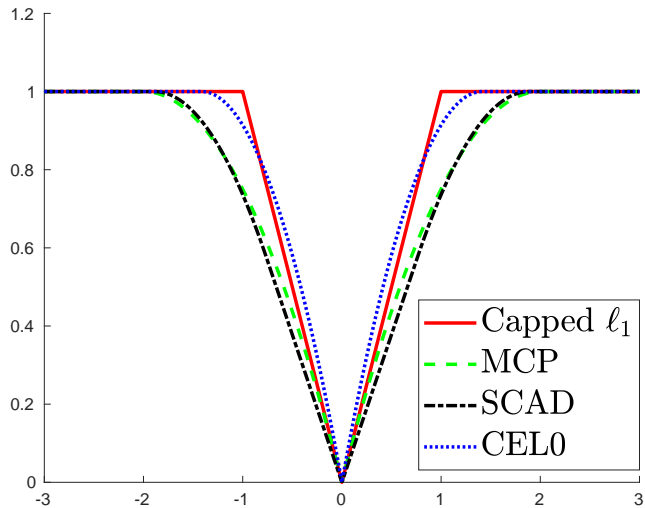}}
  \caption[Examples of continuous relaxation of $\ell_{0}$ penalization.]
          {\tabular[t]{@{}l@{}}Examples of continuous relaxation of $\ell_{0}$
            penalization \\
            ($\lambda=1$, $\gamma_{\mathrm{SCAD}}=2.5$, $\gamma_{\mathrm{MCP}}=2$,
            $\gamma_{\mathrm{CEL0}}=1$).\endtabular}
  \label{fig:reg_fct}
\end{figure}

Several examples of functions $\Psi_{\lambda}$ shown in the literature to yield good
approximations to the $\ell_{0}$ function are actually piecewise rational functions,
for which we will show in this article that exact minimization is achievable.
We list below examples of the most commonly used piecewise rational
approximations to the $\ell_{0}$ penalization that appear in several areas such as
imaging or statistics.
Figure~\ref{fig:reg_fct} displays the graph of those functions on $[-3,3]$.
\begin{itemize}
\item Capped $\ell_{p}$~\cite{Zhang_T_2010_j-mach-learn-res_analysis_mcrsr,
    Artina_M_2013_j-siam-j-optim_linearly_cnnm,Jezierska_A_2011_book_fast_stcpqcsm}:
  \begin{equation*}
    \Psi_{\lambda}(x)= \abs{x}^{p} \mathds{1}_{\left\{\abs{x} \leq \lambda\right\}}
    + \lambda^{p} \mathds{1}_{\left\{\abs{x} > \lambda\right\}} \, .
  \end{equation*}
\item Smoothly clipped absolute deviation
  (SCAD)~\cite{Fan_J_2001_j-amer-stat-asso_variable_snplop}:
  ($\gamma \in ]2,+\infty[$)
  \begin{equation*}
    \begin{aligned}
      \Psi_{\lambda}(x) = &\; \lambda\abs{x} \mathds{1}_{\left\{\abs{x} \leq \lambda\right\}}
      + \frac{(\gamma+1)\lambda^{2}}{2} \mathds{1}_{\left\{\abs{x} > \gamma\lambda\right\}}\\
      &- \frac{\lambda^2-2\gamma\lambda\abs{x}+x^{2}}{2(\gamma-1)}
      \mathds{1}_{\left\{\lambda < \abs{x} \leq \gamma\lambda\right\}} \, ,
    \end{aligned}
  \end{equation*}
\item Minimax concave penalty
  (MCP)~\cite{Zhang_C_2010_j-ann-appl-stat_nearly_uvsumcp}:
  ($\gamma \in \RR_{+}^{*}$)
  \begin{equation*}
    \Psi_{\lambda}(x)
    = \left(\lambda\abs{x}-\frac{x^{2}}{2\gamma}\right)\mathds{1}_{\left\{\abs{x} \le \gamma\lambda\right\}}
    + \frac{\gamma\lambda^{2}}{2} \mathds{1}_{\left\{\abs{x} > \gamma\lambda\right\}}  \, ,
  \end{equation*}
\item Continuous exact $\ell_{0}$
  (CEL0)~\cite{Soubies_E_2015_j-siam-j-imaging-sci_continuous_elplsrp}:
  ($\gamma \in \RR_{+}^{*}$)
  \begin{equation*}
    \Psi_{\lambda}(x) =  \lambda - \frac{\gamma^{2}}{2}{\left(\abs{x} -\frac{\sqrt{2\lambda}}{\gamma}\right)}^{2}
    \mathds{1}_{\left\{\abs{x} \leq \frac{\sqrt{2\lambda}}{\gamma}\right\}} \, .
  \end{equation*}
\end{itemize}
Although CEL0 and MCP share a similar expression for the function $\Psi_{\lambda}$, they
are quite different in their overall form $\mR_{\lambda}$ due to the choice of the
parameter $\gamma$.
In~\eqref{eq:sep_reg}, this parameter for MCP is fixed for all the samples
$x_{t}$, while for CEL0, its value is adapted to each sample.
In the above penalization, the lower the parameter $\gamma$, the tighter the
approximation to the $\ell_{0}$ penalization but the stronger also the
nonconvexity.
An important remark concerning the above examples is that, when $\Phi$ is set to
the identity, a suitable choice of the parameter $\gamma$ guarantees that the
global minimizers of the criterion $f_{\*y}+\mR_{\lambda}$ are exactly the global
minimizers of the criterion
$f_{\*y} + \lambda\ell_{0}$~\cite{Soubies_E_2017_j-siam-j-optim_unified_vecpl2l0m}.
The choice of $\gamma$ depends on the parameter $\lambda$ and the norm of the columns of
$D_{\delta}\*H$.
This behavior provides important insights and guarantees on the quality of the
above functions as penalization terms to enforce sparsity of the solutions.


\section{Rational/polynomial formulation of the problem}
\label{sec:l0_approx}

\subsection{Ubiquity of rational modeling}
\label{ssec:choice_reg}

Let us remind that the signal reconstruction problem is tackled through the
minimization of criterion $\mJ$ which has been defined
in~\eqref{eq:init_crit}.
We thus want to find
\begin{equation}
  \label{eq:orig_optim_pb}
  \mJ^{*} = \underset{\*x \in \RR^{T}}{\text{min}}
  \quad \mJ(\*x) \, .
\end{equation}
We emphasize that formulating our problem as a polynomial/rational one offers a
widely applicable framework.
First, let us show that there exists a polynomial reformulation of the $\ell_0$
criterion.
Indeed, choosing $\mR_{\lambda}=\lambda\ell_0$ in the penalization term of
criterion~\eqref{eq:init_crit}, the original problem~\eqref{eq:orig_optim_pb} can
be reformulated by using a rational function and polynomial constraints, as
follows:
\begin{equation}
  \label{eq:alt_relax}
  \begin{aligned}
    & \underset{(\*x,\bm{\xi}) \in \RR^{T}\times\RR^{T}}{\text{minimize}}
    & & \norm{\*y-D_{\delta}(\Phi(\*h \ast (\*x \odot \bm{\xi})))}^{2}
    + \lambda \sum_{t=1}^{T} \xi_{t} \\
    & \text{s.t.}
    & & (\forall t \in \nint{1,T}) \quad \xi_{t} = \xi_{t}^{2} \, ,
  \end{aligned}
\end{equation}
where the operator $\odot$ denotes the element-wise Hadamard product.
The $\xi_i$'s are introduced to formulate the $\ell_0$ penalization in a polynomial
form, while the constraints ensure that they are binary variables. 
In this formulation unfortunately, both the number of variables and the degree
of the involved polynomials are increased by a factor of two.
As a consequence, we will show in Section~\ref{sec:sdp_compl}
that~\eqref{eq:alt_relax} has a high complexity.
The method presented in this article allows us to overcome this complexity
barrier by using a different formulation of~\eqref{eq:init_crit}.

Looking closer at the relaxations of $\ell_0$ mentioned in
Section~\ref{ssec:ex_reg}, an original alternative approach consists in
considering penalization functions that are piecewise rational and can be
expressed under the general form 
\begin{equation}
  \label{eq:zbin_def}
  (\forall x \in \RR) \quad
  \Psi_{\lambda}(x) = \sum_{i=1}^{I} \zeta_{i}(x)\mathds{1}_{\{\sigma_{i-1} \leq x < \sigma_{i}\}} \, ,
\end{equation}
where ${(\zeta_{i})}_{i \in \nint{1,I}}$ are rational functions, $I$ is a nonzero integer,
and ${(\sigma_{i})}_{i \in \nint{0,I}}$ is an increasing sequence of real values.
The resulting criterion is thus a sum of rational and piecewise rational
functions:
\begin{equation}
  \label{eq:criterion}
  \begin{aligned}
    \mJ(\*x)
    &= \sum_{u=1}^{U}g_{u}(x_{\Delta(u,\alpha)-L+1},\ldots,x_{\Delta(u,\alpha)}) \\
    & \quad + \sum_{t=1}^{T}\sum_{i=1}^{I} \zeta_{i}(x_{t})\mathds{1}_{\{\sigma_{i-1} \leq x_{t} < \sigma_{i}\}} \, .
  \end{aligned}
\end{equation}


\subsection{Piecewise rational criteria}
\label{ssec:recast}

In this section, we first show how to transform the piecewise rational criterion
in Problem~\eqref{eq:orig_optim_pb} into the equivalent minimization of a sum of
rational functions under polynomial constraints.
To do so, we introduce the binary variables
${\left(z^{(i)}\right)}_{i \in \nint{1,I}}$ such that
\[
  (\forall i \in \nint{0,I}) \quad z^{(i)} = \mathds{1}_{\{\sigma_{i} \leq x\}} \, .
\]
We set $\sigma_{0}=-\infty$, $z^{(0)}=1$ and $\sigma_{I}=+\infty$, $z^{(I)}=0$ to define $\Psi_{\lambda}$ on
the whole real line $\RR$.
From the definition of ${\left(z^{(i)}\right)}_{i \in \nint{1,I}}$, we deduce that
\begin{equation}
  \label{eq:ind2bin}
  (\forall i \in \nint{0,I}) \quad
  \mathds{1}_{\{\sigma_{i-1} \leq x < \sigma_{i}\}}  = z^{(i-1)}(1-z^{(i)}).
\end{equation}
Finally, the constraint $z^{(i)} = \mathds{1}_{\{\sigma_{i} \leq x\}}$ is equivalent to two
polynomial constraints
\begin{equation}
  \label{eq:new_constr}
  z^{(i)} = \mathds{1}_{\{\sigma_{i} \leq x\}} \Longleftrightarrow \left\{
    \begin{array}{l}
      {\left(z^{(i)}\right)}^{2} - z^{(i)} = 0 \\
      \left(z^{(i)}-\frac{1}{2}\right)\left(x-\sigma_{i}\right) \geq 0  \, . \\
    \end{array}
  \right.
\end{equation}

Indeed, the polynomial equality constraint enforces $z^{(i)}$ to be a binary
variable while the polynomial inequality constraint ensures that it takes the
same values as $\mathds{1}_{\{\sigma_{i} \leq x\}}$ for every $x$ in $\RR$.
Therefore, by substituting~\eqref{eq:new_constr} for~\eqref{eq:criterion},
Problem~\eqref{eq:orig_optim_pb} reads as the minimization of a sum of rational
functions depending on $\*x$ and vectors
$\*z = {\left(\*z^{(i)}\right)}_{i\in\nint{0,I}}$ under polynomial constraints,
namely
\begin{equation}
  \begin{aligned}
    \label{eq:recast_pb}
    & \underset{\substack{\*x,\*z}}{\text{\rm minimize}}
    & & \sum_{u=1}^{U}g_{u}(x_{\Delta(u,\alpha)-L+1},\ldots,x_{\Delta(u,\alpha)}) \\
    & & & + \sum_{t=1}^{T}\sum_{i=1}^{I} \zeta_{i}(x_{t})z_{t}^{(i-1)}(1-z_{t}^{(i)}) \\
    & \text{s.t.}
    & & (\forall (i,t) \in \nint{0,I} \times \nint{1,T}) \quad
    \begin{cases}
      {\left(z_{t}^{(i)}\right)}^{2} - z_{t}^{(i)} = 0 \\
      \left(z_{t}^{(i)}-\frac12\right)\left(x_{t}-\sigma_{i+1}\right) \geq 0
      \, .
    \end{cases}
  \end{aligned}
\end{equation}

More generally, this reformulation can be applied to the minimization of any
piecewise rational function.
For instance, a piecewise rational fit function $f_{\*y}$ could also be chosen.


\subsection{Symmetry of regularizers}
\label{ssec:sym_reg}

All the piecewise rational approximations to the $\ell_{0}$ penalty listed in
Section~\ref{ssec:ex_reg} are even functions.
This symmetry property is expressed here by an absolute value on the input
variable $x_{t}$ in the expressions of function $\Psi_{\lambda}$.
This absolute value is handled in our framework by adding an additional variable
$r_{t}$ for each $x_{t}$ and adding the two constraints
\begin{equation}
  \label{eq:abs_constr}
  \left\{
    \begin{aligned}
      r_{t}^{2} &= x_{t}^{2} \\
      r_{t}  &\geq 0 \, .
    \end{aligned}
  \right.
\end{equation}
This symmetry is important to decrease the number $I$ of variables $\*z$
involved in~\eqref{eq:recast_pb} and therefore to reduce the overall complexity
of the final problem to be solved, as will be explained in
Section~\ref{sec:sdp_compl}.
Indeed, it can divide by two the number of pieces in $\Psi_{\lambda}$, leading to only
$I/2$ pieces instead of $I$.
Taking the example of the MCP penalization, instead of having the four intervals,
$]-\infty,-\gamma\lambda[$, $[-\gamma\lambda,0[$, $[0,\gamma\lambda[$, and
$[\gamma\lambda,+\infty[$, we have only the two intervals,
$[0,\gamma\lambda[$ and $[\gamma\lambda,+\infty[$.
Using symmetry results in adding one variable $r_{t}$ and $I/2$ variables
${\left(z_{t}^{(i)}\right)}_{i \in \nint{1,I/2}}$ for each $x_{t}$ as well as
$2+I/2$ polynomial constraints corresponding to
constraints~\eqref{eq:new_constr} and~\eqref{eq:abs_constr}.
This has to be compared with the direct formulation where we introduce $I$
variables ${\left(z_{t}^{(i)}\right)}_{i \in \nint{1,I}}$ with $I$ polynomial
constraints.
Note that in our analysis of Section~\ref{sec:sdp_compl}, we omit the equality
constraints that force ${\left(z_{t}^{(i)}\right)}_{i \in \nint{1,I/2}}$ to be
binary variables since substitution will be performed for those constraints in
Section~\ref{ssec:substitution}.


\section{Solving the optimization problem}
\label{sec:solve}

This section is concerned with the resolution of Problem~\eqref{eq:recast_pb}
presented in a progressive manner.
After a brief review of techniques from polynomial and rational optimization in
Section~\ref{ssec:relax_ratio}, we apply the latter to our signal processing
context.
In Sections~\ref{ssec:relax_sum_ration} and~\ref{ssec:relax_ourpb}, we explicitly
show how the structure of our problem allows us to reduce the dimensions of the
final convex relaxation.
Our analysis reveals that signal processing problems are computationally
tractable when using sparsity patterns and subsampling.


\subsection{Minimizing a rational function}
\label{ssec:relax_ratio}

A sum of rational functions can be written as a single rational function by
reduction to a common denominator.
A first step to handle~\eqref{eq:recast_pb} is hence to consider the
minimization of a single rational function.
In this section, we simplify our notation to explain the framework used to
solve~\eqref{eq:recast_pb} and we focus on the generic problem of finding
\begin{equation}
  \label{eq:ex_ratio}
  \mJ^{*} = \min_{\*x \in \mK} \; \frac{p(\*x)}{q(\*x)} \, ,
\end{equation}
where $p$ and $q$ are polynomials in $T$ variables and $\mK \subset \RR^T$ is
the feasible set.
Section~\ref{ssec:relax_ourpb} will get back to Problem~\eqref{eq:recast_pb}.


\subsubsection{Condition on the feasible set}

In~\eqref{eq:ex_ratio}, $\mK$ is a basic subset of $\RR^{T}$ defined by polynomial
inequalities as
\begin{equation}
  \label{eq:def_fes-set}
  \mK = \{\*x \in \RR^{T} \mid (\forall j \in \nint{1,J}) \quad s_{j}(\*x) \geq 0\}
  \, ,
\end{equation}
where, for every $j \in \nint{1,J}$, $s_j\colon \RR^T \to \RR$.
As we often work with bounded signals, we make the mild assumption that $\mK$
contains $T$ polynomial constraints of the form
\[
  (\forall t \in \nint{1,T}) \quad x_{t}^{2} \leq B^{2} \, ,
\]
where $B$ is a positive constant.
Since $\mK$ is a closed set in a finite dimensional space, the above boundedness
condition ensures that $\mK$ is a compact set.
To simplify the notation, we write those constraints into a vector form as
\[
  (\*B-\*x) \odot (\*x+\*B) \geq \*0 \, ,
\]
where $\*B$ and $\*0$ are the vectors composed solely of $B$ and $0$,
respectively.


\subsubsection{Reformulation as a moment problem}

As shown in~\cite[Proposition 5.20]{Lasserre_J_2009_book_moments_ppta},
Problem~\eqref{eq:ex_ratio} is equivalent to find
\begin{equation}
  \label{eq:meas_eq}
  \begin{aligned}
    \inf_{\mu \in \mM_{+}(\mK)} & \int_{\mK}p(\*x)\mu(\diff\*x) \\
    \text{s.t.} & \int_{\mK}q(\*x)\mu(\diff\*x) = 1 \, ,
  \end{aligned}
\end{equation}
where $\mM_{+}(\mK)$ denotes the set of positive finite measures supported on
$\mK$.
The equivalence between Problems~\eqref{eq:ex_ratio} and~\eqref{eq:meas_eq} relies
on the possibility to link any optimal point $\*x_{*}$ of~\eqref{eq:ex_ratio} to
a Dirac measure $\delta(\*x_{*})/q(\*x_{*})$ solution to~\eqref{eq:meas_eq}.
The main idea here is to embed the original problem in a higher dimensional
space in order to linearize it.
At first glance,~\eqref{eq:meas_eq} looks more intricate than~\eqref{eq:ex_ratio}
since we need to minimize over an infinite-dimensional set of measures supported
by $\mK$ instead of minimizing on $\mK$ itself.
However, the objective function and the constraint are now linear in the new
optimized variable $\mu$.
Furthermore, by defining $\*x^{\bm{\alpha}}=x_{1}^{\alpha_{1}} \dots x_{T}^{\alpha_{T}}$, notice that
\begin{equation}
  \label{eq:reform_obj}
  \int_{\mK}p(\*x)\mu(\diff \*x)
  = \int_{\mK}\sum_{\bm{\alpha} \in \NN^{T}} p_{\bm{\alpha}}\*x^{\bm{\alpha}}\mu(\diff \*x)
  = \sum_{\bm{\alpha} \in \NN^{T}} p_{\bm{\alpha}} v_{\bm{\alpha}} \, ,
\end{equation}
where $v_{\bm{\alpha}} = \int_{\mK}\*x^{\bm{\alpha}}\mu(\diff \*x)$ denotes the
moment of order $\bm{\alpha}$ of the measure $\mu$.
For convenience, we will use infinite vectors to write sums such as the
rightmost member of~\eqref{eq:reform_obj}.
We define the infinite vector $\widetilde{\*p}={(p_{\bm{\alpha}})}_{\bm{\alpha} \in \NN^{T}}$
and the infinite moment vector $\widetilde{\*v}={(v_{\bm{\alpha}})}_{\bm{\alpha} \in \NN^{T}}$.
Since $\widetilde{\*p}$ has a finite number of nonzero elements, the sum
in~\eqref{eq:reform_obj} is well defined and can be written
$\widetilde{\*p}^{\top}\widetilde{\*v}$.

Since $\mK$ is a compact set, the measure $\mu$ is uniquely defined by its moments
and thus we can reformulate Problem~\eqref{eq:meas_eq} as
\begin{equation}
  \label{eq:mom_eq}
  \begin{aligned}
    & \inf_{\widetilde{\*v}\in \RR^{\NN^{T}}}
    & & \widetilde{\*p}^{\top}\widetilde{\*v} \\
    & \text{s.t.} & & \widetilde{\*q}^{\top}\widetilde{\*v} = 1 \\
    & & & \widetilde{\*v} \in \mD(\mK)
  \end{aligned} \, ,
\end{equation}
where $\widetilde{\*q}$ is defined similarly to $\widetilde{\*p}$ as the
infinite vector extensions of $\*q$ obtained by zero padding and $\mD(\mK)$ is
the cone of moments of positive measures supported on $\mK$.
Our objective now is to replace this difficult conic constraint by simpler
constraints.
We introduce two tools, respectively, the moment matrix $\*M(\widetilde{\*v})$
associated to the moment vector $\widetilde{\*v}$ and the localizing matrix
$\*M^{s}(\widetilde{\*v})$ associated to $\widetilde{\*v}$ with respect to a given
polynomial $s$.
Those matrices are infinite-dimensional and are both defined through their
entries as follows
\begin{align*}
  \left(\forall (\bm{\alpha},\bm{\beta}) \in \NN^{T} \times \NN^{T}\right) \quad
  M_{\bm{\alpha},\bm{\beta}}(\*v) &= v_{\bm{\alpha}+\bm{\beta}} \\
  \left(\forall (\bm{\alpha},\bm{\beta}) \in \NN^{T} \times \NN^{T}\right) \quad
  M_{\bm{\alpha},\bm{\beta}}^{s}(\*v) &=
  \sum_{\bm{\gamma} \in \NN^{T}} s_{\bm{\gamma}}v_{\bm{\alpha}+\bm{\beta}+\bm{\gamma}}
  \, .
\end{align*}
We define such infinite-dimensional matrices to be positive semi-definite if
all their finite-dimensional principal submatrices are positive semi-definite.
Since $\mK$ is compact, Putinar's
theorem~\cite[Proposition 3.1]{Henrion_D_2013_misc_optimization_lmipsc} states
that $\widetilde{\*v}$ has a representing measure in $\mM_{+}(\mK)$ if and only
if the corresponding moment matrix $\*M(\widetilde{\*v})$ and localizing
matrices ${(\*M^{s_{j}}(\widetilde{\*v}))}_{j\in\nint{1,J}}$ are positive
semi-definite.
The positive semi-definiteness of the moment and localizing matrices guarantee
that $\widetilde{\*v}$ represents a positive measure and ensures that its
support is $\mK$.


\subsubsection{Converging hierarchy of SDP problems}

To solve numerically Problem~\eqref{eq:mom_eq}, we replace the conic constraint
with semidefinite constraints and then truncate the moment vector
$\widetilde{\*v}$, as well as its associated moment and localizing matrices, up
to a degree $2k$ for a given integer $k$.
This yields a hierarchy of convex SDP problems, known as Lasserre's
hierarchy~\cite{Lasserre_J_2001_j-siam-j-optim_global_oppm}.
For a given relaxation order $k$, the SDP relaxation to be solved reads:
\begin{equation}
  \label{eq:sdp_eq}
  \begin{aligned}
    & \mJ^{*}_{k} = &&\inf_{\*v \in \RR^{m}}
    & & \sum_{\bm{\alpha} \in \NN_{2k}^{T}} p_{\bm{\alpha}}v_{\bm{\alpha}} \\
    & && \text{s.t.}
    & & \sum_{\bm{\alpha} \in \NN_{2k}^{T}} q_{\bm{\alpha}}v_{\bm{\alpha}} = 1 \\
    & && & & \*M_{k}(\*v) \in \Sym_{+}^{n_{0}} \\
    & && & & (\forall j \in \nint{1,J}) \quad
    \*M^{s_{j}}_{k-\mathrm{d}_{s_{j}}}(\*v) \in \Sym_{+}^{n_{j}} \, ,
  \end{aligned} 
\end{equation}
where ${(\mathrm{d}_{s_{j}})}_{j\in\nint{1,J}}$ are defined
in~\eqref{eq:defdemideg}.
The cardinality of $\NN_{2k}^{T}$ is $\binom{T+2k}{2k}$ and thus $\*v$ is a
vector containing $m=\binom{T+2k}{2k}$ moments.
Furthermore, the truncated moment matrix $\*M_{k}$ is the principal submatrix of
the infinite-dimensional moment matrix $\*M$ that has dimension
$n_{0} \times n_{0}$ with $n_{0}=\binom{T+k}{k}$.
Thereby, $\*M_{k}$ is indexed by $(\bm{\alpha},\bm{\beta})$ in
$\NN_{k}^{T} \times \NN_{k}^{T}$ and contains all the moments up to degree $2k$.
Similarly, the truncated localizing matrices are the submatrix of their
infinite-dimensional counterparts that have the dimensions $n_{j} \times n_{j}$
with $n_{j}=\binom{T+k-\mathrm{d}_{s_{j}}}{k-\mathrm{d}_{s_{j}}}$.

Problem~\eqref{eq:sdp_eq} is an SDP problem in its dual form with linear
equality constraints.
Indeed, aggregating the moment and localizing matrices into a single symmetric
block diagonal matrix before separating it into a sum along the elements of
$\*v$, we obtain
\begin{equation}
  \label{eq:sdp_pb}
  \begin{aligned}
    & \mJ^{*}_{k} = &&\underset{\*v \in \RR^{m}}{\text{minimize}}
    & & \*p^{\top} \*v \\
    & && \text{s.t.} \quad
    & & \*C-\sum_{i=1}^{m} v_{i} \*A_{i} \in \Sym_{+}^{n} \\
    & && & & \*a - \*G^{\top} \*v = \*0  \, , 
  \end{aligned}
\end{equation}
where $\*C$ and ${(\*A_{i})}_{i \in \nint{1,m}}$ are symmetric matrices, $\*a$
is a vector of $ \RR^{\ell}$, and $\*G$ is a matrix of $\RR^{m \times \ell}$.
The dimension $n$ is thus given by $n=\sum_{j=0}^{J} n_{j}$.
Notice that, in this section,~\eqref{eq:sdp_eq} has only one linear constraint
thus $\ell=1$, $\*a=1$, and $\*G = \*q$.
However, in next sections, more linear constraints will be involved, so that we
prefer to employ this matrix-vector notation here.

Solving each SDP problem yields a lower bound $\mJ^{*}_{k}$ on the optimal value
$\mJ^{*}$ of the criterion $\mJ$.
Furthermore, the higher the order $k$, the tighter the bound $\mJ^{*}_{k}$ but
the higher also the dimensions of the SDP problem.
In our context where the sought signal is bounded, ${(\mJ^{*}_{k})}_{k \in \NN}$
is an increasing convergent sequence whose limit is
$\mJ^{*}$~\cite{Lasserre_J_2001_j-siam-j-optim_global_oppm}.
Moreover, the hierarchy has finite convergence generically, i.e.\ convergence
happens at a finite relaxation
order generically~\cite{Nie_J_2013_j-math-prog_optimality_cfclh}.
Finally, an exact global solution $\hat{\*x}$ of~\eqref{eq:ex_ratio} can be
extracted from the solution of an SDP problem~\eqref{eq:sdp_eq} indexed by a
relaxation order at which convergence has
occurred~\cite{Henrion_D_2005_book_detecting_goesg}.

Note that the relaxation order $k$ should be chosen such that
\[
  k \geq \max\left\{\mathrm{d}_{p},\mathrm{d}_{q},
  \max_{j \in \nint{1,J}}\mathrm{d}_{s_{j}}\right\}
  \, .
\]
This is a necessary condition which ensures that $2k$ is greater than the
maximum degree of $p$, $q$ and all the ${(s_{j})}_{j \in \nint{1,J}}$, and
prevents truncation of the latter polynomials.
There is no a priori known sufficient relaxation order to ensure the convergence
of the hierarchy.
However, once the SDP relaxation is solved, there exists a sufficient condition
that guarantees the convergence.
Namely, if the moment matrices $\*M_{k}$ and $\*M_{k-1}$ have the same rank,
then convergence has occurred~\cite{Lasserre_J_2009_book_moments_ppta}.

We remark that the dimensions $n$ and $m$ of the SDP problem grows respectively
as $T^{k}$ and $T^{2k}$ when $T$ is large, hence exponentially in the degree of the
involved polynomials.


\subsection{Problem structure emerging from a sum of rational
  functions}
\label{ssec:relax_sum_ration}

Although constituting the theoretical foundation of our work, the approach
presented in Section~\ref{ssec:relax_ratio} is computationally inefficient and
requires further improvements that we now explain.
Indeed, reducing a sum of rational functions using a common denominator often
yields a rational function with very high degree, which then requires a high
relaxation order $k$ in the hierarchy of SDP problems.
As a consequence, the obtained SDP problems are too high-dimensional to be
solvable in a reasonable time using state-of-the-art solvers.
However, a more ingenious method is to use the structure induced by the sum
to yield a block SDP
problem~\cite{Castella_M_2019_j-ieee-trans-sig-proc_rational_onral0p}.
There are two types of structure to consider in our problem: first we deal with
a sum of rational functions instead of a single one, and then each of those
functions has only a few subset of variables as input.
Those two structures are sometimes referred to as sparse problem and sparse
polynomials~\cite{Waki_H_2006_j-siam-j-optim_sums_ssprpopss,
  Bugarin_F_2015_j-math-prog-comput_minimizing_smrf}.
However, in order to prevent confusion with the sparsity of the original signal
$\bxini$, we will not use this terminology.
To illustrate it, let us turn our attention on finding
\begin{equation}
  \label{eq:ex_sum_ratio}
  \mJ^{*} = \min_{\*x \in \mK}
  \sum_{i=1}^{\tilde{I}}\frac{p_{i}(\*x_{E_{i}})}{q_{i}(\*x_{E_{i}})} \, ,
\end{equation}
where $p_{i}$ and $q_{i}$ are polynomials in $T_{i}$ variables and $\mK$ is a compact
subset of $\RR^{T}$ having the form~\eqref{eq:def_fes-set}.
The vector $\*x_{E_{i}}$ denotes the subvector of $\*x$ composed of the elements
indexed by the set $E_{i}$, the set $E_{i}$ being a subset of $\nint{1,T}$ of
cardinality $T_{i}$.
We further assume that the polynomials in~\eqref{eq:ex_sum_ratio} involve only a
few variables, i.e.\
\begin{equation}
  \label{eq:sparse_assumpt}
  \left(\forall i \in \nint{1,\tilde{I}}\right) \quad T_{i} \ll T \, .
\end{equation}


\subsubsection{Exploiting the sum of rational functions structure}
\label{sssec:exploit_struc}

Instead of introducing a single measure on all the variables, we now introduce
a measure $\mu_{i}$ for each rational function $p_{i}/q_{i}$ of the sum.
However, coupling between variables from different measures appears when the
sets ${(E_{i})}_{i \in \nint{1,\tilde{I}}}$ intersect.
We therefore need to add moment equality constraints to ensure that the
overlapping moments of two measures are identical.
Moreover, we need some restrictions on how variables can overlap several
measures in order to keep the problem consistent.
The required condition is that the sets ${(E_{i})}_{i \in \nint{1,\tilde{I}}}$ verify
the so-called running intersection
property~\cite{Lasserre_J_2009_book_moments_ppta,
  Bugarin_F_2015_j-math-prog-comput_minimizing_smrf} which is stated as
\[
  (\forall i \in \nint{2,\tilde{I}}, \exists \tilde{\jmath} \in \nint{1,i-1}) \quad
  E_{i} \bigcap \left(\bigcup_{j=1}^{i-1}E_{j}\right) \subseteq E_{\tilde{\jmath}} \, .
\]
Together with the compactness of $\mK$, it guarantees that
Problem~\eqref{eq:ex_sum_ratio} is equivalent to
\begin{equation}
  \label{eq:sum_mom_relax}
  \begin{aligned}
    \inf_{\bm{\mu} \in \Xi}
    & \sum_{i=1}^{\tilde{I}}\int_{\mK_{i}}p_{i}(\*x_{E_{i}})\mu_{i}(\diff\*x_{E_{i}}) \\
    \text{s.t.}\; & (\forall i \in \nint{1,\tilde{I}}) \,
    \int_{\mK_{i}}q_{i}(\*x_{E_{i}})\mu_{i}(\diff\*x_{E_{i}}) = 1 \\
    & (\forall (i,j) \in \nint{1,\tilde{I}}\times\nint{1,\tilde{I}})
    (\forall \bm{\gamma} \in \NN^{\Card(E_{i,j})}) \\
    & \int_{\mathcal{K}_{i}}q_{i}(\*x_{E_{i}})\*x_{E_{i,j}}^{\bm{\gamma}}\mu_{i}(\diff\*x_{E_{i}})
    = \int_{\mathcal{K}_{j}}q_{j}(\*x_{E_{j}})\*x_{E_{i,j}}^{\bm{\gamma}}\mu_{j}(\diff\*x_{E_{j}})
    \, ,
  \end{aligned}
\end{equation}
where $E_{i,j} = E_{i} \bigcap E_{j}$ and $\bm{\mu}={(\mu_{i})}_{i \in \nint{1,\tilde{I}}}$ is the
new optimization variable belonging to the product
$\Xi=\bigtimes_{i\in\nint{1,\tilde{I}}}\mM_{+}(\mK_{i})$.
The sets ${(\mK_{i})}_{i \in \nint{1,\tilde{I}}}$ are subsets of $\RR^{T_{i}}$ defined
by the subsets of polynomials in variables $\*x_{E_{i}}$ defining $\mK$.
The last equality constraints in~\eqref{eq:sum_mom_relax} enforce equality
between the marginal distributions of $q_{j}\mu_{j}$ and $q_{i}\mu_{i}$ along
$\*x_{E_{ij}}$.
In other words, those constraints ensure the equality of overlapping moments
between the different measures.


\subsubsection{Block structure in the SDP hierarchy}

As in Section~\ref{ssec:relax_ratio}, we now use Putinar's theorem to replace
each measure by its moment vector at the cost of additional semi-definite
constraints.
We then truncate the moment vectors as well as the moment and
localizing matrices, before stacking them.
As a result, the moment vector
$\*v = {[\*v_{1}^{\top},\ldots,\*v_{\tilde{I}}^{\top}]}^{\top}$ is a stack of the moment
vectors of each measure $\mu_{i}$.
Similarly, the moment matrix $\*M_{k}(\*v)
=\Diag\left(\*M_{1,k}(\*v_{1}),\dots,\*M_{\tilde{I},k}(\*v_{\tilde{I}})\right)$
and the localizing matrices
$\*M^{s_{j}}_{k-\mathrm{d}_{s_{j}}}(\*v)
= \Diag\Big(\*M^{s_{j}}_{1,k-\mathrm{d}_{s_{j}}}(\*v_{1}),\dots
,\*M^{s_{j}}_{\tilde{I},k-\mathrm{d}_{s_{j}}}(\*v_{\tilde{I}})\Big)$ have a
block diagonal structure where each diagonal block corresponds respectively to
the moment or localizing matrix of one of the measures $\mu_{i}$.
This leads to the following SDP problem:
\begin{equation}
  \label{eq:sdp_eq_sum_ration}
  \begin{aligned}
    \mJ^{*}_{k} = &\inf_{\*v \in \RR^{m}} {\*p}^{\top}\*v \\
    \text{s.t.} & \; (\forall i \in \nint{1,\tilde{I}}) \quad
    {\*q_{i}}^{\top}\*v_{i} = 1 \\
    & \*M_{k}(\*v) \in \Sym_{+}^{n_{0}} \\
    & (\forall j \in \nint{1,J}) \quad
    \*M^{s_{j}}_{k-\mathrm{d}_{s_{i}}}(\*v) \in \Sym_{+}^{n_{j}} \\
    & \*F\*v = \*0
    \, ,
  \end{aligned}
\end{equation}
where
\[
  \begin{aligned}
    m    &= \sum_{i=1}^{I}\binom{T_{i}+2k}{2k} \, , \;
    n_{0} = \sum_{i=1}^{I}\binom{T_{i}+k}{k} \, ,  \;
    n_{j} =
    \sum_{i=1}^{I}\binom{T_{i}+k-\mathrm{d}_{s_{j}}}{k-\mathrm{d}_{s_{j}}} \, , \\
    \*p &= {[\*p_{1}^{\top}, \ldots, \*p_{\tilde{I}}^{\top}]}^{\top} \, , \\
  \end{aligned}
\]
and $\*F$ is a matrix in $\RR^{\ell \times m}$ representing the linear
constraints linking the ${(\*v_{i})}_{i \in \nint{1,\tilde{I}}}$ together and coming
from the constraints between the projections in~\eqref{eq:sum_mom_relax}.
Similarly to Section~\ref{ssec:relax_ratio},~\eqref{eq:sdp_eq_sum_ration} can be
finally expressed in the canonical form~\eqref{eq:sdp_pb}.

There are two main differences with the situation discussed in
Section~\ref{ssec:relax_ratio}:
\begin{itemize}
\item Instead of having a single measure on all the variables, we obtain several
  measures on different smaller subsets of variables.
  The SDP optimization variable $\*v$ is now a vector built by stacking the
  different truncated moment vectors ${(\*v_{i})}_{i \in \nint{1,\tilde{I}}}$ of each
  measure.
  As a consequence, the moment and localizing matrices have a block diagonal
  structure, each block corresponding to a measure, or equivalently to a term in
  the sum of Problem~\eqref{eq:ex_sum_ratio}.
  Thanks to Assumption~\eqref{eq:sparse_assumpt}, the size of the blocks in the
  moment matrix, equal to $\binom{T_{i}+2k}{2k}$, is much smaller than the size
  $\binom{T+2k}{2k}$ obtained in Section~\ref{ssec:relax_ratio}.
  Especially when $k$ increases, the difference in size becomes even more
  significant.
  The block structure can then be efficiently exploited by SDP solvers to
  decrease the computational time.
\item Extra moment constraints due to the coupling between variables arise.
  Although those constraints may be numerous, they are linear equality
  constraints in the SDP problem; their impact on the computational time of the
  SDP solver is minor.
\end{itemize}


\subsection{Minimizing our criterion}
\label{ssec:relax_ourpb}

We now apply the method of Section~\ref{ssec:relax_sum_ration} to
solve~\eqref{eq:recast_pb}.
We have to handle a sum of $\tilde{I} = U+T$ terms.
We hence introduce a measure for each term, i.e.\ $U$ measures
${(\mu_u)}_{u \in \nint{1,U}}$ for the rational functions ${(g_u)}_{u \in \nint{1,U}}$
and $T$ measures ${(\nu_{t})}_{t \in \nint{1,T}}$ for the rational functions in the
reformulated penalization.
The measures ${(\mu_u)}_{u \in \nint{1,U}}$ are measures on at most $L$ variables
$x_{\Delta(u,\alpha)-L+1},\ldots,x_{\Delta(u,\alpha)}$ while the measures
${(\nu_{t})}_{t \in \nint{1,T}}$ are measures on $I+1$ scalar variables,
corresponding to $x_{t}$ and $\*z_{t}$.


\subsubsection{Feasible set for our reformulated problem}

In Problem~\eqref{eq:recast_pb}, the sets ${(\mK_{i})}_{i \in \nint{1,U+T}}$ are
defined by the bound constraints and the polynomial constraints arising from the
reformulation of Section~\ref{ssec:recast}.
Namely, the sets ${(\mK_{i})}_{i \in \nint{1,U}}$ are defined by
\begin{equation}
  \label{eq:supp_cons_1}
  (\forall i \in \nint{1,U}) \quad
  (\*B-\*x_{E_{i}})\odot(\*x_{E_{i}}+\*B) \geq \*0 \, ,
\end{equation}
while the sets ${(\mK_{i})}_{i \in \nint{U+1,U+T}}$ are defined by
\begin{equation}
  \label{eq:supp_cons_2}
  \begin{aligned}
    (\forall i  \in \nint{U+1,U+T}) \quad
    &(B-x_{i-U})(x_{i-U}+B) \geq 0 \\
    &(\forall j \in \nint{1,I}) \quad
    {(z_{i-U}^{(j)})}^{2} - z_{i-U}^{(j)} = 0 \\
    &(\forall j \in \nint{1,I}) \quad
    \left(z_{i-U}^{(j)}-\frac{1}{2}\right)\left(x_{i-U}-\sigma_{j+1}\right) \geq 0
    \, .
  \end{aligned}
\end{equation}
Note that, since we introduce a measure for each rational function in the sum,
we have to cope with more than $T$ bound constraints.
Indeed several measures are defined on identical variables and we need to
introduce bound constraints for each of those measures.
We then perform the relaxation~\eqref{eq:sdp_eq_sum_ration} to generate a
hierarchy of SDP problems.


\subsubsection{Coupling and linear equality constraints}
\label{sssec:coupling}

We observe that two kinds of coupling as discussed in
Section~\ref{sssec:exploit_struc} appear: one between the different measures
${(\mu_u)}_{u \in \nint{1,U}}$ and one between the measures ${(\mu_u)}_{u \in \nint{1,U}}$
and the measures ${(\nu_{t})}_{t \in \nint{1,T}}$.
By definition of the convolution matrix $\*H$ in Section~\ref{sec:motivation},
the sets ${(E_{i})}_{i \in \nint{1,\tilde{I}}}$ satisfy the running intersection
property.
Furthermore, among the extra moment equality constraints to acknowledge
coupling, we remark that many of them between moments of the measures
${(\mu_u)}_{u \in \nint{1,U}}$ are redundant.
Let us take a simple example to illustrate this fact.

\paragraph*{Example}
Assume that we want to minimize over the variable
$\*x={(x_{t})}_{t \in \nint{1,5}}$, a sum of three rational functions which has
the following form
\[
  \frac{p_{1}(x_{1},x_{2},x_{3})}{q_{1}(x_{1},x_{2},x_{3})} + \frac{p_{2}(x_{2},x_{3},x_{4})}{q_{2}(x_{2},x_{3},x_{4})}
  + \frac{p_{3}(x_{3},x_{4},x_{5})}{q_{3}(x_{3},x_{4},x_{5})}
  \, .
\]
Following the method developed in Section~\ref{ssec:relax_sum_ration}, we
introduce three measures $\mu_{1}$ $\mu_{2}$ and $\mu_{3}$, one for each term of the sum.
We thus need the following equality constraints between moments, for every
$(\alpha,\beta)$ in $\NN^{2}$,
\[
  \begin{aligned}
    \int q_{1}(x_{1},x_{2},x_{3})x_{2}^{\alpha}x_{3}^{\beta} \mu_{1}(\diff x_{1}, \diff x_{2}, \diff x_{3})
    &= \int q_{2}(x_{2},x_{3},x_{4})x_{2}^{\alpha}x_{3}^{\beta} \mu_{2}(\diff x_{2}, \diff x_{3}, \diff x_{4}) \\
    \int q_{2}(x_{2},x_{3},x_{4})x_{3}^{\alpha}x_{4}^{\beta} \mu_{2}(\diff x_{2}, \diff x_{3}, \diff x_{4})
    &= \int q_{3}(x_{3},x_{4},x_{5})x_{3}^{\alpha}x_{4}^{\beta} \mu_{3}(\diff x_{3}, \diff x_{4},\diff x_{5}) \\
    \int q_{1}(x_{1},x_{2},x_{3})x_{3}^{\alpha} \mu_{1}(\diff x_{1}, \diff x_{2}, \diff x_{3})
    &= \int q_{3}(x_{3},x_{4},x_{5})x_{3}^{\alpha} \mu_{3}(\diff x_{3}, \diff x_{4},\diff x_{5})
    \, .
  \end{aligned}
\]
We observe that the variable $x_{3}$ appears in each term of the sum and thus
also in moments of each measure.
In particular, we notice that the last constraint is redundant with
the first two ones when $\beta=0$.
It is thus sufficient to consider only moment equality constraints on
consecutive measures ${(\mu_u)}_{u \in \nint{1,U}}$ in~\eqref{eq:recast_pb}.
We can thus drastically reduce the number of moment linear equality constraints.


\section{Complexity of the relaxation}
\label{sec:sdp_compl}

Current state-of-the-art SDP solvers use interior points methods which are known
to be very efficient for small and medium scale problems.
On the other hand, their running time becomes prohibitive for large scale
problems.
This is a major drawback of the relaxation of polynomial optimization problems
into SDP problems.
Nevertheless, Sections~\ref{ssec:relax_sum_ration} and~\ref{ssec:relax_ourpb} used
the structure of the problem to yield a structured SDP problem.
In this section, we derive the complexity of this SDP problem and show that it
is computationally solvable in a fair amount of time.

The complexity of an SDP problem under the form~\eqref{eq:sdp_pb} is expressed
as a quadruple of integers $(n,m,m_{\mathrm{s}},\ell)$.
The integer $m$ denotes the size of the vector of optimized variables, $n$ is
the size of the semi-definite inequality constraint, $\ell$ is the number of linear
equality constraints, and $m_{\mathrm{s}}$ is the number of block matrices
involved in the semi-definite constraint.
Note that $n$ is related to $m_{\mathrm{s}}$ since it is the sum of the size of
each block.
The above quadruple therefore does not fully characterize the structure of an
SDP\@.
For example, having one huge block and nine tiny ones is not equivalent in terms
of complexity to having ten medium blocks.
However, knowing $n$ and $m_{\mathrm{s}}$ is usually enough to get a good
evaluation of the complexity of the problem.

The bottleneck for current SDP solver is mainly the dimension of both $n$ and
$m$.
This section gives an asymptotic estimation for $n$ and $m$ according to the
parameters of our initial model~\eqref{eq:model} and the relaxation order $k$.
A more detailed derivation for the expression of $(n,m,m_{\mathrm{s}},\ell)$ is
presented in Appendix~\ref{sec:detail_complex_comp}.
We show first that subsampling and sparsity allow to decrease the latter and
make the numerical resolution of the associated SDP problems tractable.
Then, we introduce tools and tricks that allow us to decrease further the
dimension of the SDP problems to be solved, so reducing the computational
time of our method.


\subsection{Consequence of the subsampling on the dimensions of the SDP problem}
\label{ssec:ben_subsampl}

For a given relaxation order $k$, when the number of samples $T$ goes to
infinity and $L \gg 1$ (i.e.\ we lose the band structure of $\*H$), the size
of the SDP problem asymptotically becomes of the order
(see Appendix~\ref{sec:detail_complex_comp})
\begin{equation}
  \label{eq:asymp_compl}
  m = \mathcal{O}(UL^{2k}+T) \quad , \quad
  n = \mathcal{O}(UL^{k}+T) \, .
\end{equation}
We note that both sizes $n$ and $m$ grow exponentially with $k$ and blow up
quickly.
In particular, $m$ grows faster than $n$.
However, we will see that the SDP hierarchy often converges quickly in practice,
that is ${\left(\mJ_{k}^{*}\right)}_{k \in \NN}$ converge to $\mJ^{*}$ for a
relaxation order $k$ of $2$, $3$, or $4$.
From our analysis, we observe that the main bottleneck of our method is the
number of variables per measure and the order of relaxation.
While the number of variables in measures ${(\nu_{t})}_{t \in \nint{1,T}}$ is fixed
to $I+1$, the total number of variables in measures ${(\mu_u)}_{u \in \nint{1,U}}$ is
$L$ and~\eqref{eq:asymp_compl} shows that $m$ and $n$ rise quickly with $L$.

Although the subsampling reduces the quality of the reconstruction by
eliminating some information on the signal, as a side effect in our context, it
allows the size of the SDP relaxation to be reduced.
As shown by~\eqref{eq:Udelta}, decimation decreases $U$, which plays a prominent
role in the complexity parameters $(n,m,m_{\mathrm{s}},\ell)$ of the SDP problem.
Table~\ref{table:sdp_relax_size} compares the size of SDP relaxations for the
SCAD penalization without decimation ($D_{\infty}$) and with $D_{2}$ (resp. $D_{4}$)
decimation.
As discussed above, the dimensions $n$ and $m$ increase quickly with the
relaxation order $k$ and the length of the filter $L$.
Note that because of the approximation made in Section~\ref{ssec:nb_blocks},
stating that measures ${(\mu_u)}_{u \in \nint{1,U}}$ are on $L$ variables, the SDP
dimension presented here are slightly overestimated.

\begin{table}[htbp]
  \caption{Dimension of the relaxation of the SCAD penalization for different
    decimations}
  \begin{center}
    \setlength\tabcolsep{3.5pt}
    \resizebox{0.95\textwidth}{!}{\begin{tabular}{ccccccccccccccc}
      \cline{4-15}
      \multicolumn{3}{c}{} & \multicolumn{3}{c}{$m$} & \multicolumn{3}{c}{$n$}
      & \multicolumn{3}{c}{$m_{\mathrm{s}}$} & \multicolumn{3}{c}{$\ell$} \\
      \midrule
      $T$ & $L$ & $k$ & $D_{\infty}$ & $D_{4}$ & $D_{2}$ & $D_{\infty}$ & $D_{4}$ & $D_{2}$
                      & $D_{\infty}$ & $D_{4}$ & $D_{2}$ & $D_{\infty}$ & $D_{4}$ & $D_{2}$ \\
      \midrule
       50 &  3  &  3  &  8400 &  7476 &  6300 &  7000 &  6450 &  5750
                      &   600 &   556 &   500 &  1035 &   735 &   420 \\
      \midrule
      100 &  3  &  3  & 16800 & 14784 & 12600 & 14000 & 12800 & 11500
                      &  1200 &  1104 &  1000 &  2085 &  1755 &   845 \\
      \midrule
       50 &  4  &  3  & 14700 & 12390 &  9450 &  9250 &  8205 &  6875
                      &   650 &   595 &   425 &  2015 &  1355 &   660 \\
      \midrule
      100 &  4  &  3  & 29400 & 24360 & 18900 & 18500 & 16220 & 13750
                      &  1300 &  1180 &  1050 &  4065 &  3405 &  1335 \\
      \midrule
      100 &  5  &  3  & 54600 & 43512 & 31500 & 25600 & 21236 & 17050
                      &  1400 &  1256 &  1100 &  7530 &  6375 &  2315 \\
      \midrule
       50 &  3  &  4  & 16500 & 14685 & 12375 & 13500 & 12455 & 11125
                      &   600 &   556 &   500 &  1772 &  1184 &   568 \\
      \midrule
      100 &  3  &  4  & 33000 & 29040 & 24750 & 27000 & 24720 & 22250
                      &  1200 &  1104 &  1000 &  3572 &  2102 &  1143 \\
      \midrule
      100 &  4  &  4  & 66000 & 54120 & 41250 & 38500 & 33460 & 28000
                      &  1300 &  1180 &  1050 &  9116 &  7268 &  2172 \\
      \bottomrule
    \end{tabular}}
  \end{center}
  ~\label{table:sdp_relax_size}
\end{table}


\subsection{Polynomial equality constraints and substitution}
\label{ssec:substitution}

For a given measure, equality constraints involving monic monomials in the
definition of the support set $\mK_{i}$ can be substituted.
The constraint is then used to reduce the number of moments in the vector of
moments.
We clarify this process here through the example of the SCAD penalization.
Substitution is carried out automatically by some
software~\cite{Henrion_D_2009_j-optim-meth-softw_gloptipoly_mosdp}, but has not
been clearly documented.

Let us focus our attention on the measure $\nu_{t}$, depending on the three
variables $x_{t}$, $z_{t}^{(1)}$, and $z_{t}^{(2)}$, as well as on the
associated truncated vector $\*v_{t}$ of moment up to degree $2$.
Using the equality constraints in~\eqref{eq:supp_cons_2}, we substitute the
related monomial in $\*v_{t}$.
The moments associated with monomials ${\left(z_{t}^{(1)}\right)}^{2}$ and
${\left(z_{t}^{(2)}\right)}^{2}$ are thus the same as the ones associated with
$z_{t}^{(1)}$ and $z_{t}^{(2)}$.
Therefore, the moment vector $\*v_{t}$ has a dimension reduced by two.
When $\*v_{t}$ contains moments up to degree $2k$, substitution reduces the
number of moments from $\binom{3+2k}{2k}$ to $8k$.

In the general case, for a given relaxation order $k$, $\*v_{t}$ contains only
$2k(k+1)$ moments after substitution which is much fewer than the original
$\binom{1+I+2k}{2k}$ moments.
Substitution significantly decreases the values of $n$, $m$ and
$m_{\mathrm{s}}$ which have a major impact on the computational cost of SDP
solvers.
However, it does not impact the number of linear constraints $\ell$.


\subsection{Linear versus quadratic polynomial constraints}
\label{ssec:lin_vs_quad}

Our method to solve rational optimization problem is valid only if the
constraint set $\mK$ is compact.
We therefore set a bound $B$ on the sought signal.
The bound constraints can be expressed in two ways:
\begin{itemize}
\item first, as two linear vector constraints
  \[
    \begin{aligned}
      \*B - \*x &\geq \*0  \\
      \*x + \*B &\geq \*0 \, ,
    \end{aligned}
  \]
\item or as a single quadratic vector constraint
  \[
     (\*B-\*x)\odot(\*x+\*B) \geq \*0 \, .
  \]
\end{itemize}
Following Appendix~\ref{ssec:nb_blocks}, using two linear inequality constraints
per variable introduces $2(U+L)$ localizing matrices and consequently $2(U+L)$
blocks in our SDP problems while using a quadratic inequality constraint only
adds $U+L$ blocks.
Moreover, linear and quadratic constraints yield blocks of identical size.
Indeed, the size of a localizing matrix $\*M_{k}^{s}$ corresponding to a polynomial
$s$ in $\omega$ variables is given by $\binom{\omega+k-\mathrm{d}_{s}}{k-\mathrm{d}_{s}}$
and here $\mathrm{d}_{s}=1$ for both linear and quadratic constraints.
Therefore, formulating the bound constraints as quadratic constraints reduces by
a factor two the number of blocks associated to such bounds.


\subsection{Using a sign oracle}
\label{ssec:oracle}

For real-valued signals $\bxini$, convergence is observed at orders $k$ for
which building and solving the corresponding SDP problems is highly demanding in
terms of computation and memory storage.
Conversely, when $\bxini$ is a positive signal, we
observed~\cite{Castella_M_2019_j-ieee-trans-sig-proc_rational_onral0p} convergence
at a lower order $k$.
This suggests a method yielding similar results for real-valued signals using an
oracle.
Instead of
\begin{equation*}
  (\forall \*x \in \RR^{T}) \quad
  \mJ(\*x) =
  \frac{1}{2}\norm{\*y-D_{\delta}(\Phi(\*H\*x))}^{2} + \sum_{t=1}^{T} \Psi_{\lambda} (x_{t})
  \, ,
\end{equation*}
we minimize
\begin{equation*}
  (\forall \*x \in \RR_{+}^{T}) \quad
  \widetilde{\mJ}(\*x) =
  \frac{1}{2}\norm{\*y-D_{\delta}(\Phi(\widetilde{\*H}\*x))}^{2} + \sum_{t=1}^{T} \Psi_{\lambda} (x_{t})
  \, ,
\end{equation*}
where $\widetilde{\*H}=\*H\Diag({\bm{\epsilon}})$, $\bm{\epsilon}\in\{-1,1\}^{T}$ is the
sign vector of $\bxini$ provided by the oracle, and $\Diag({\bm{\epsilon}})$ is a
diagonal matrix with binary elements $\bm{\epsilon}$.
We build our oracle by solving a standard least absolute shrinkage and
selection operator (LASSO)
problem~\cite{Tibshirani_R_1996_j-r-stat-soc-b-stat-meth_regression_ssl}, i.e.\
$\Psi_{\lambda}=\lambda |.|$.
The availability of an oracle allows us to restrict the minimization
of~\eqref{eq:criterion} to positive valued signals thanks to the new
convolution matrix $\widetilde{\*H}$.
Our oracle decreases significantly the computational time in two ways:
\begin{itemize}
\item Since the convergence of the SDP hierarchy occurs for smaller order $k$,
  the dimensions of the SDP problem to solve are much lower according to
  Section~\ref{sec:sdp_compl}.
\item Moreover, since we optimize now on positive variables, we do not need to
  use the additional variables ${(r_{t})}_{t \in \nint{1,T}}$ introduced in
  Section~\ref{ssec:sym_reg} to account for symmetries and the presence of
  absolute values.
  This results in smaller vectors of moments, hence a lower dimensional SDP
  problem.
\end{itemize}
An exact solution is thus retrieved by solving an SDP problem of fair dimension.
Finally, the computational cost of our oracle is low since we solve a
LASSO using a forward-backward algorithm.
It typically takes less than a second which is negligible compared to the
computational time of our method as shown in Section~\ref{sec:simulations} while
providing accurate oracle on the sign of the initial signal.


\section{Numerical simulations and results}
\label{sec:simulations}

\subsection{Experimental set-up}

To show the efficiency of our framework, we apply it to the reconstruction
of a sparse signal subject to nonlinear distortion and subsampling.
We use a piecewise relaxation of $\ell_{0}$ to promote sparsity as detailed in
Section~\ref{ssec:ex_reg}.
We perform simulations on 50 test cases where the initial sparse signal
$\bxini$ has length $T=100$ with $10$ non-zero values.
Those values are drawn randomly according to a uniform distribution on
$[-1,-0.1]\cup[0.1,1]$.
The position of the non-zero values are also drawn randomly according to uniform
distribution on $\nint{1,T}$.
The length $L$ of the filter is set to $3$ and its coefficients are the
normalized $L$-th row of Pascal's triangle.
This kind of filter is useful to model enlargement due to measurement from
sensors for example.
We choose the following saturation function for the nonlinear distortion $\phi$
\[
  (\forall t \in \mathbb{R}) \quad \phi(t) = \frac{t}{\chi + |t|} \, ,
\]
where $\chi$ is set to $0.3$.
Finally, we perform the relaxation into SDP for relaxation orders $2$, $3$, and
$4$.
We use GloptiPoly~\cite{Henrion_D_2009_j-optim-meth-softw_gloptipoly_mosdp}
to relax rational problems into SDP problems which are then solved with the
solver SDPT3~\cite{Toh_K_1999_j-optim-meth-softw_sdpt3_mspspv13}.
All the simulations have been run on a standard computer with an Intel Xeon CPU
running at 3.7~GHz and 32 GB of RAM allocated to the process.


\subsection{Example of a rational relaxation: SCAD}
\label{ssec:ex_ratio_refom}

To clarify the reformulation of Section~\ref{ssec:recast}, we demonstrate it on
the regularizers given in Section~\ref{ssec:ex_reg}.
Taking advantage of symmetry as explained in Section~\ref{ssec:sym_reg}, SCAD has
three pieces and thus requires to introduce variables $z_{t}^{(1)}$
and $z_{t}^{(2)}$ leading to
\begin{equation}
  \label{eq:ex_scad}
  \begin{aligned}
    & \underset{\substack{\*x,\*z}}{\text{\rm minimize}}
    & & f_{\*y}(\*x)
    + \sum_{t=1}^{T} (1-z_{t}^{(1)})\lambda \abs{x_{t}}  + z_{t}^{(2)} \frac{(\gamma+1)\lambda^{2}}{2} \\
    & & & - z_{t}^{(1)}(1-z_{t}^{(2)})\frac{\lambda^{2}-2\gamma\lambda\abs{x_{t}}+x_{t}^{2}}{2(\gamma-1)} \\
    & \text{s.t.}
    & & (\forall (i,t) \in \{1,2\} \times \nint{1,T}) \quad
    {\left(z_{t}^{(i)}\right)}^{2} - z_{t}^{(i)} = 0 \\
    & & & (\forall t \in \nint{1,T}) \quad
    \left(z_{t}^{(1)}-\frac{1}{2}\right)\left(\abs{x_{t}}-\gamma\lambda\right) \geq 0 \\
    & & & (\forall t \in \nint{1,T}) \quad
    \left(z_{t}^{(2)}-\frac{1}{2}\right)\left(\abs{x_{t}}-\lambda\right) \geq 0 \, .
  \end{aligned}
\end{equation}

A similar approach applies to Capped $\ell_{p}$, MCP, and CEL0 penalties;
the details are omitted for conciseness.
Although we use SCAD penalization in all the subsequent simulations, similar
results can be obtained with Capped $\ell_{p}$, MCP, and CEL0.
Nonetheless, SCAD is more demanding in terms of computation since it has more
rational pieces.
It consequently provides a worst case scenario for the computational time
compared with the other penalizations.
The parameter $\gamma$ for SCAD is set to $2.1$ in order to approximate $\ell_{0}$
closely.
The value of the parameter $\lambda$ was determined empirically and set to $0.15$.


\subsection{Acceleration of convergence with the sign oracle}
\label{ssec:oracle_simu}

In this section, we want to show how the oracle impacts the convergence of the
SDP hierarchy.
We first consider the use of a sign oracle in a linear model, i.e.\ the case
when $\phi=\Id$.
We then delve into the more challenging case of a nonlinear model.
Decimation is set to $D_{4}$ in this section.
The oracle is build on solving a LASSO problem by using a forward-backward
algorithm as described in Section~\ref{ssec:oracle}.


\subsubsection{Linear case}
\label{sssec:oracle_lin}

\begin{figure}[!t]
  \centering
  \captionsetup{justification=centering}
  \subfloat[With oracle: $k = 2$ (top), $3$ (middle) $4$ (bottom)]
           {\includegraphics[width=0.80\linewidth]{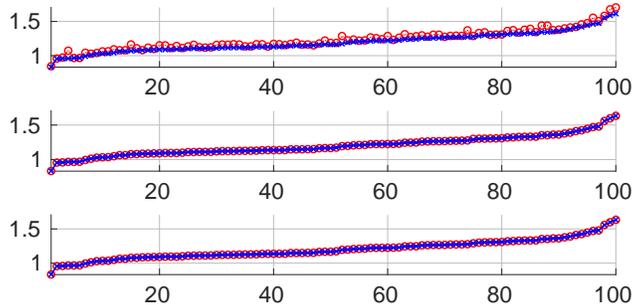}}
           
  \subfloat[Without oracle: $k = 2$ (top), $3$ (middle) $4$ (bottom)]
           {\includegraphics[width=0.80\linewidth]{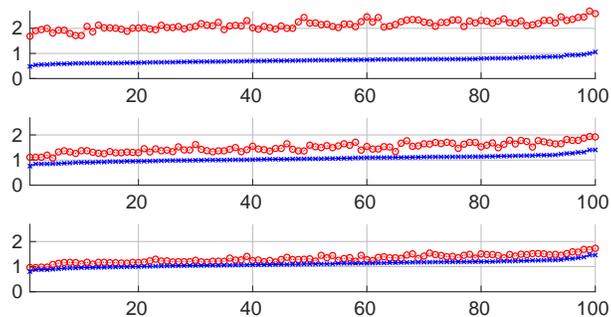}}
                    
  \caption{Comparison between the lower bound $\mJ^{*}_{k}$ and the value of the
    criterion $\mJ(\hat{\*x}_{k})$ for 100 tests (Linear Case).}
  ~\label{fig:crit_vs_binf_lin}
\end{figure}

Solving each SDP problem in the hierarchy provides both a lower bound
$\mJ^{*}_{k}$, which is the value of the objective function of the SDP at
optimality, and an approximate minimizer $\hat{\*x}_{k}$, which is extracted
from a minimizer of the SDP problem.
We compare here the value of the criterion at $\hat{\*x}_{k}$ with
$\mJ^{*}_{k}$.
Since increasing the relaxation order $k$ yields larger lower bounds and smaller
criterion values, we know that the convergence of the hierarchy happens when
$\mJ(\hat{\*x}_{k})$ and $\mJ^{*}_{k}$ are equal.
Figure~\ref{fig:crit_vs_binf_lin} compares those two values respectively in the
cases with oracle and without the use of our oracle on 100 test cases.
From top to bottom, the two figures are drawn for relaxation orders $k=2$,
$k=3$, and $k=4$.
Criterion values are represented in red while lower bounds are represented in
blue.
Each point of the $x$-axis represents the values for a single test case.
For the sake of clarity, the values are ordered according to the value of the
lower bound.
We observe that, without oracle, the convergence is slow and still not reached
in general at order $k=4$.
On the other hand, when we use our oracle, convergence appears quickly, i.e.\
$k=3$ in most of the test cases.

\subsubsection{Nonlinear case}
\label{sssec:oracle_nonlin}

\begin{figure}[!ht]
  \centerline{\includegraphics[width=0.8\linewidth]{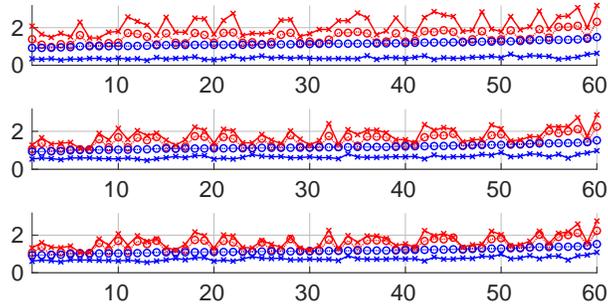}}
  \caption{Comparison between the lower bound $\mJ^{*}_{k}$ and the value of the
    criterion $\mJ(\hat{\*x}_{k})$ when the oracle is used.
    Plain line with cross dots: no oracle used,
    dashed line with circle dots: use of our oracle (Nonlinear case).}
  ~\label{fig:crit_vs_binf_nonlin}
\end{figure}

Figure~\ref{fig:crit_vs_binf_nonlin} is similar to
Figure~\ref{fig:crit_vs_binf_lin} but in the context of a nonlinear model.
The continuous line with cross dots represents the cases without the use of an
oracle while the dashed line with circle dots represents the cases with our sign
oracle.
We observe here that even with a sign oracle, the convergence of the hierarchy
does not occur for low values of $k$ due to the nonlinearity.
However, we can notice that the gap between the lower bound and the criterion
value at the $\hat{\*x}_{k}$ is greatly reduced when we use our oracle.


\subsection{Reconstruction of sparse signals}

\subsubsection{Global optimality}
\label{sssec:glob_opti}

In this section, we want to demonstrate the quality of the minimizers of
Problem~\eqref{eq:ex_scad} returned by various methods.
Note that we do not use the oracle here.
We use the decimation operator $D_{4}$ but similar results hold for the other
operators.
We compare our method to a forward-backward (FB) algorithm applied directly to
the criterion $\mJ = f_{\*y}+\mR_{\lambda}$ where the gradient step is first performed
on the data fitting term and a proximal step is then performed on the
penalization.
Hence the criterion to minimize is the same for both methods.
We initialize the FB algorithm first with the null vector and denote by
$\*x_{\mathrm{FB}0}$ the resulting solution.
Then we perform a warm start of the FB algorithm using the solution obtained
from our method as an initializer.
The resulting estimate is denoted by $\*x_{\mathrm{FB}1}$.

\begin{figure}[!ht]
  \centerline{\includegraphics[width=0.8\linewidth]{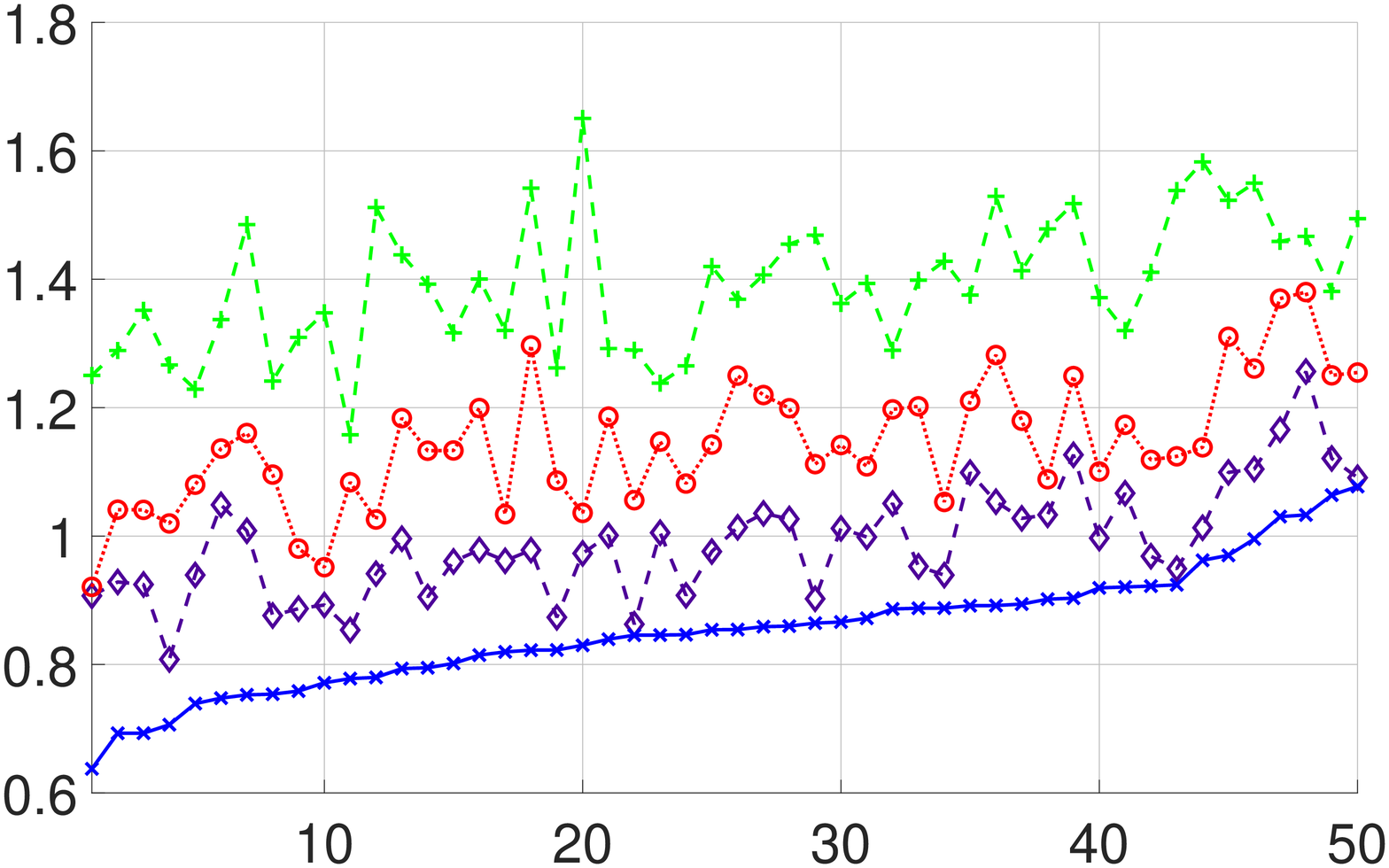}}
  \caption{Comparison between the different values of the criterion for the
    minimizers returned by the different methods.
    In red $\mJ(\hat{\*x}_{4})$, in blue $\mJ^{*}_{4}$, in green
    $\mJ(\*x_{\mathrm{FB}0})$, and in purple $\mJ(\*x_{\mathrm{FB}1})$.}
  ~\label{fig:comp_crit_warmstart}
\end{figure}

In Figure~\ref{fig:comp_crit_warmstart}, we compare the value of the criterion
$\mJ$ at $\*x_{\mathrm{FB}0}$ and $\*x_{\mathrm{FB}1}$ with the solution
returned by our method for a relaxation order $k=4$.
The solid blue curve with cross dots represents the values of the lower bound
$\mJ^{*}_{4}$, the pointed red curve with circle dots represents $\mJ(\hat{\*x}_{4})$,
the dashed green curve with plus dots represents $\mJ(\*x_{\mathrm{FB}0})$, and
the dashed purple curve with plus dots represents $\mJ(\*x_{\mathrm{FB}1})$.

Since the criterion $\mJ$ is highly nonconvex, the forward-backward algorithm
gets stuck in local minimizers.
Indeed, changing the initialization point changes the output of the algorithm.
We can observe it on Figure~\ref{fig:comp_crit_warmstart} where the green and
purple curves are not superposed.
Moreover, similarly to Section~\ref{sssec:oracle_nonlin}, we observe that the
convergence in the hierarchy has not occurred at order $4$ since the blue and
red curves are not superimposed.
As a consequence, $\hat{\*x}_{4}$ is not a global minimizer of $\mJ$ but only an
approximation of it.
A solution to improve the quality of the minimizer is to use the solution
$\hat{\*x}_{4}$ as a warm start of the FB algorithm as shown by the purple curve.

\subsubsection{Quality of signal reconstruction}

We now look at the quality of the signal reconstruction in terms of mean square
error: our method is compared with several other ones to illustrate its interest
for faithful recovery of the original signal $\bxini$.
In addition to the FB algorithm presented in Section~\ref{sssec:glob_opti}, we
compare our method with the oracle to iLASSO, a LASSO approach modified to
handle the nonlinearity of the model.
It consists first on applying the LASSO using a linearization of the nonlinear
operator $\phi$.
Namely, it solves
\begin{equation*}
  \underset{\*x \in \RR^{T}}{\argmin}
  \norm{\*y-D_\alpha(\mL_{\phi}(\*h \ast \*x))}^{2} + \lambda_{\mathrm{LASSO}} \norm{\*x}_{1}
  \, ,
\end{equation*}
where $\mL_{\phi}$ is a linearization of $\phi$ and $\lambda_{\mathrm{LASSO}}$ is a parameter
set empirically to $0.1$.
Note that for our choice of $\phi$, $\mL_{\phi}=\chi^{-1}$.
We subsequently apply a modified iterative hard thresholding (IHT) that handles
the nonlinearity.
Namely, we apply the FB algorithm to find
\[
  \underset{\*x \in \RR^{T}}{\argmin}
  \norm{\*y-D_\alpha(\Phi(\*h \ast \*x))}^{2} + \lambda_{\mathrm{IHT}} \ell_{0}(\*x)
  \, ,
\]
where we perform a gradient step on the data fidelity component and a proximal
step on the penalization $\lambda_{\mathrm{IHT}}\ell_{0}$.
This method provides better results than the FB algorithm presented in
Section~\ref{sssec:glob_opti}.
We also compare our method to the Iteratively Reweighted $\ell_{1}$ algorithm
(IRL1)~\cite{Candes_E_2008_j-four-anal-appl_enhancing_srlm} applied to
\begin{equation*}
  \label{eq:lin_phi}
  \underset{\*x \in \RR^{T}}{\argmin}
  \norm{\*y-D_\alpha(\mL_{\phi}(\*h \ast \*x))}^{2} + \mR_{\lambda_{\mathrm{IRL1}}}(\*x)
  \, ,
\end{equation*}
where $\mR_{\lambda}$ is the SCAD regularization.
Both IRL1 and FB algorithms are initialized with the null vector.

\begin{figure}[!ht]
  \centerline{\includegraphics[width=\linewidth]{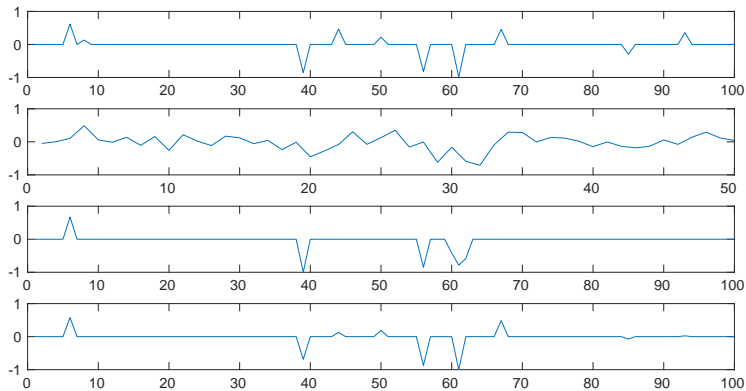}}
  \caption{Comparison between iLASSO and our method for signal reconstruction
    under nonlinear transformation and subsampling.
    From top to bottom: the original signal $\bxini$, the observed signal $\*y$,
    and respectively the signal reconstructed with iLASSO
    $\*x_{\mathrm{iLASSO}}$ and with our method $\hat{\*x}_{3}$.}
  ~\label{fig:ex_sig_recons}
\end{figure}

Figure~\ref{fig:ex_sig_recons} illustrates the different signals for a single
realization using $D_{2}$ decimation.
From top to bottom, we display the original signal $\bxini$, the subsampled
observed signal $\*y$, the signal reconstructed respectively with iLASSO
$\*x_{\mathrm{iLASSO}}$, and the signal reconstructed using our method
$\hat{\*x}_{3}$ at the relaxation order $k=3$.
We do not display the signal reconstructed with FB and IRL1 since those
algorithms are not well suited for solving~\eqref{eq:criterion} and thus
provide poor quality reconstruction.
We first notice that iLASSO misses many peaks and also detects a peak that does
not exist in the original signal while our method detects almost all peaks.
One could argue the threshold coefficient $\lambda_{\mathrm{IHT}}$ in iLASSO is too
high but, when we decrease it, small artifacts appear.
In contrast, our method detects almost all peaks and do not leave any artifact.
We observe that some peaks do not have the same amplitude as the ones in the
original signal.
This is due to subsampling.
Indeed, if a peak is located on an even index, it will be eliminated by the
subsampling.
However, the convolution with $h$, that represents the physical limitation of
sensors in our example, allows us still to recover the peak since it gets
enlarged to odd neighboring.
Even though, we lose information about the amplitude of this peak.

\begin{figure}[!t]
  \centering
  \captionsetup{justification=centering}
  \subfloat[$D_{\infty}$ decimation]
  {\includegraphics[width=0.5\linewidth]{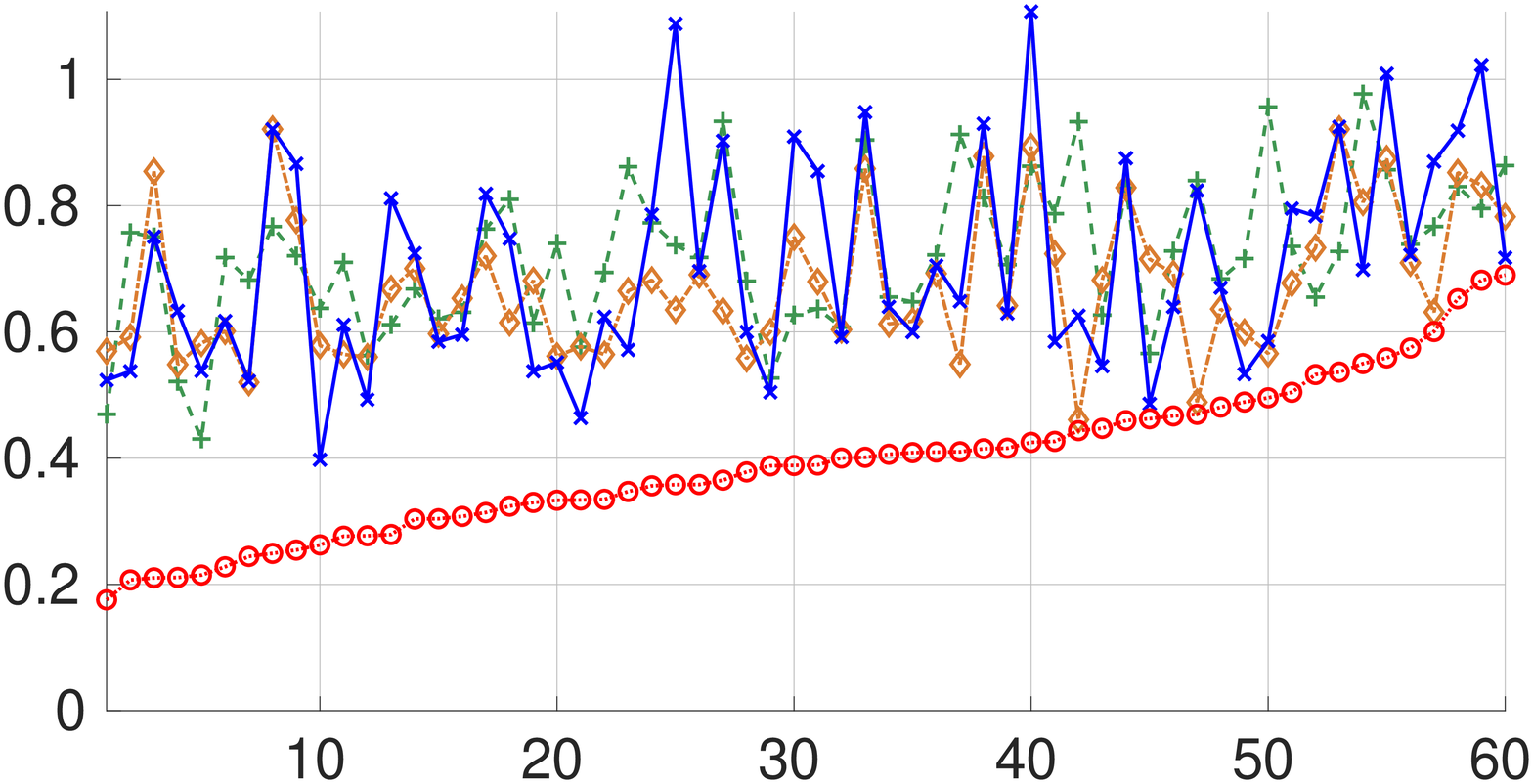}}
  \subfloat[$D_{4}$ decimation]
  {\includegraphics[width=0.5\linewidth]{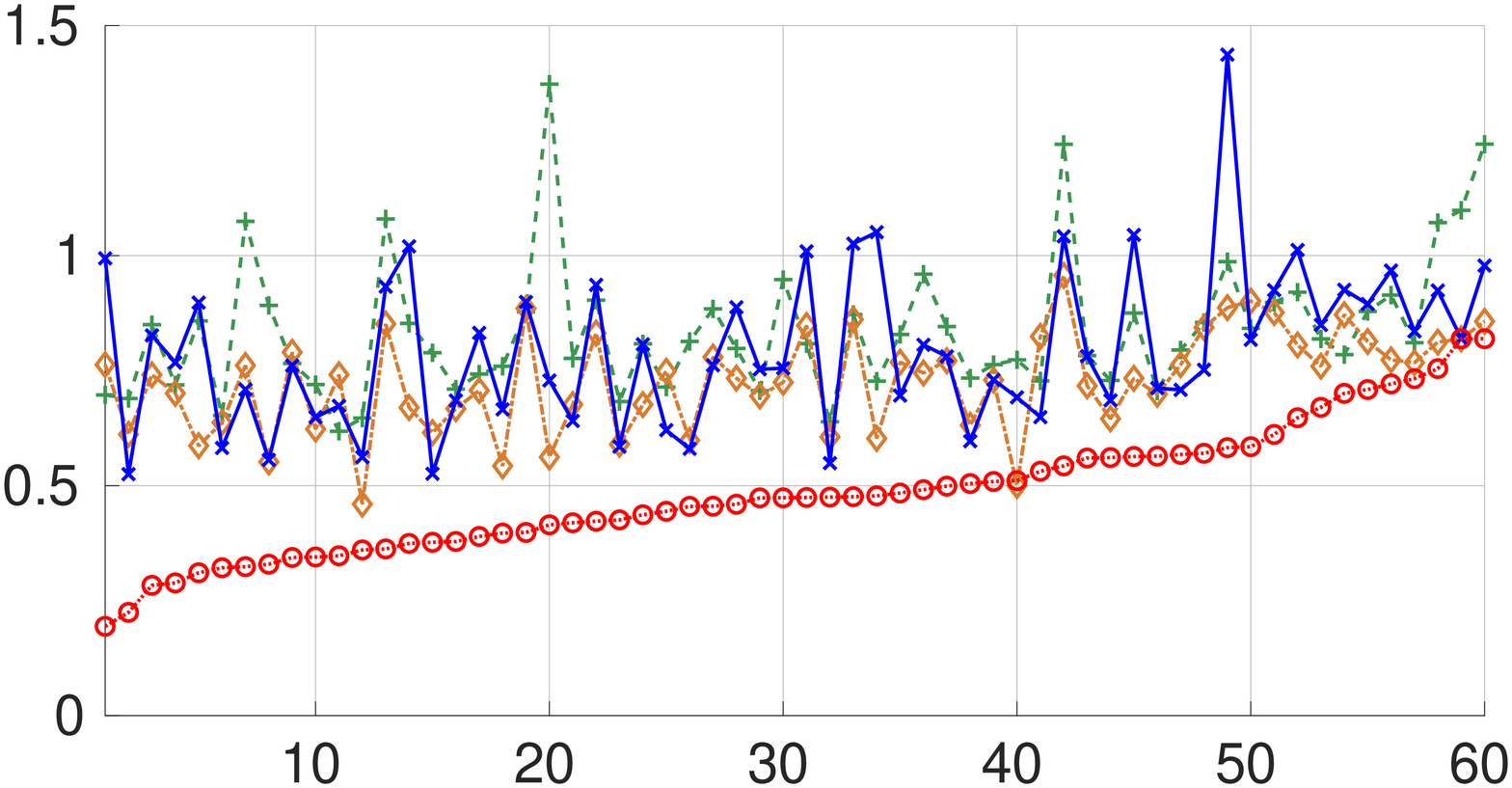}}
           
  \subfloat[$D_{2}$ decimation]{
  \raisebox{-.5\height}{\includegraphics[width=0.5\linewidth]{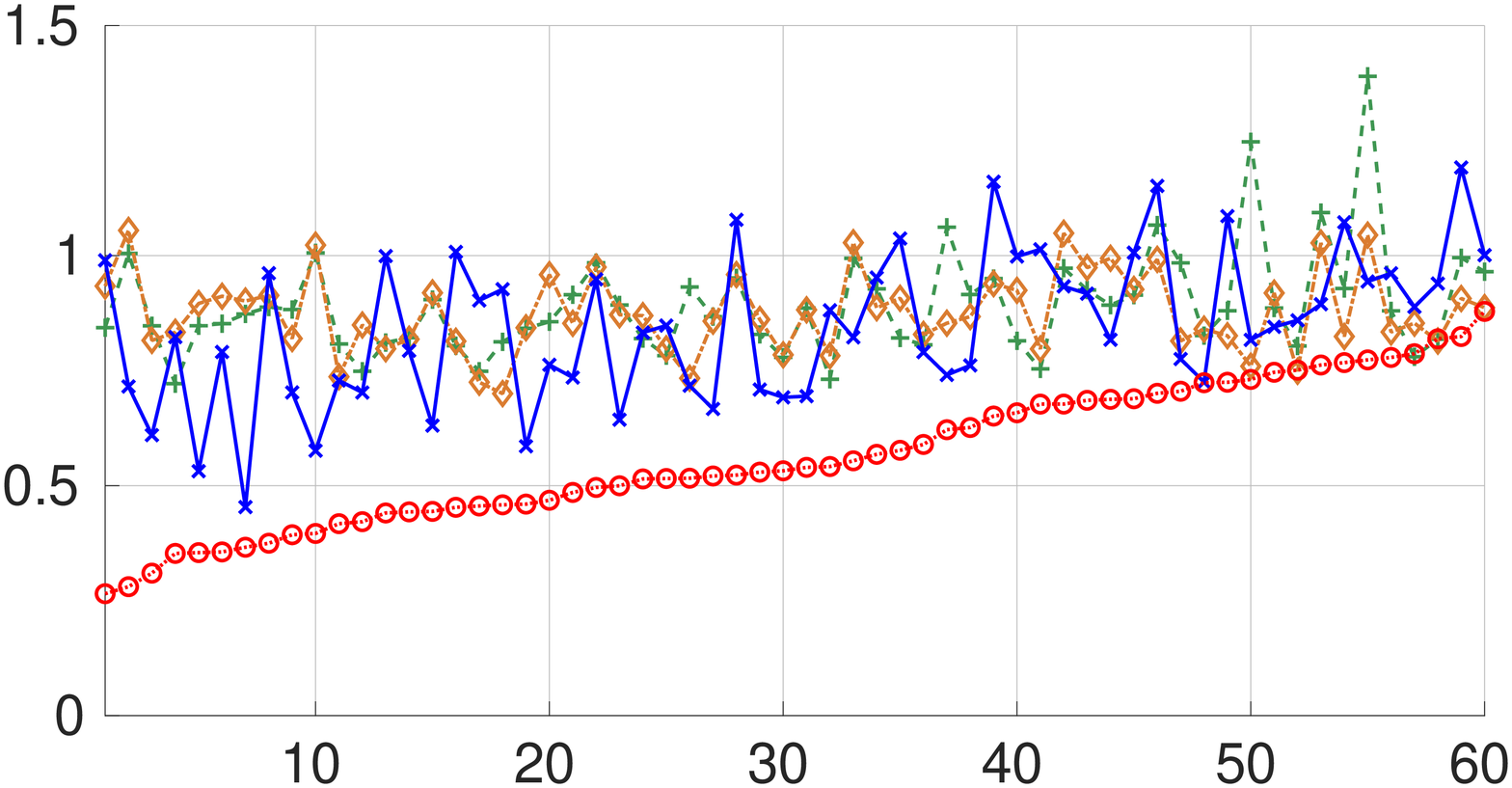}}}
  \resizebox{0.48\textwidth}{!}{
    \begin{tabular}{lccc}
      Average values     & $D_{\infty}$ & $D_{4}$ & $D_{2}$ \\
      \toprule
      $\*x_{\mathrm{FB}}$    & 0.72 & 0.84 & 0.89 \\
      $\*x_{\mathrm{IRL1}}$  & 0.67 & 0.73 & 0.87 \\
      $\*x_{\mathrm{iLASSO}}$ & 0.71 & 0.80 & 0.84 \\
      $\hat{\*x}_{4}$        & 0.39 & 0.48 & 0.56 \\
      \bottomrule
    \end{tabular}
  }
                    
  \caption{Mean square error between the estimated signal and the original
    signal $\bxini$.
    In dashed green: $\*x_{\mathrm{FB}}$, in dashed orange: $\*x_{\mathrm{IRL1}}$,
    in blue: $\*x_{\mathrm{iLASSO}}$, and in dotted red: $\hat{\*x}_{4}$.
    Average values are shown in the table.}
  ~\label{fig:mse}
\end{figure}

Figures~\ref{fig:mse} shows the mean square error
$\norm{\bxini-\*x}/\norm{\bxini}$ for $D_{\infty}$, $D_{4}$ and $D_{2}$ decimation
between the original signal $\bxini$ and: in green $\*x_{\mathrm{FB}}$, in
orange $\*x_{\mathrm{IRL1}}$, in blue $\*x_{\mathrm{iLASSO}}$, and in red
$\hat{\*x}_{4}$.
Those confirm the good reconstruction result shown in the specific example of
Figure~\ref{fig:ex_sig_recons}.

Finally, Table~\ref{table:comput_time} shows the average computational times for
different decimation operators and relaxation orders.
As we expected, the better performance of our method comes at the expense of a
higher computational cost than iLASSO, which takes less than 1 second.
\begin{table}[htbp]
  \caption{Computation time of our method (in seconds)}
  \begin{center}
    \setlength\tabcolsep{3.5pt}
    \resizebox{0.63\textwidth}{!}{\begin{tabular}{cccccccc}
        & \multicolumn{3}{c}{Without oracle} & & \multicolumn{3}{c}{With oracle} \\
        & $D_{\infty}$ & $D_{4}$ & $D_{2}$ & & $D_{\infty}$ & $D_{4}$ & $D_{2}$ \\
        \toprule
        $k=2$ &    41 &    35 &   29 &&    38 &    31 &   25 \\
        $k=3$ &   162 &   121 &   87 &&   144 &   106 &   74 \\
        $k=4$ & 29991 & 14575 & 5801 && 24362 & 11062 & 4084 \\
        \bottomrule
    \end{tabular}}
  \end{center}
  ~\label{table:comput_time}
\end{table}


\subsubsection{Handling higher-dimensional signal}

\begin{figure}[!ht]
  \centerline{\includegraphics[width=\linewidth]{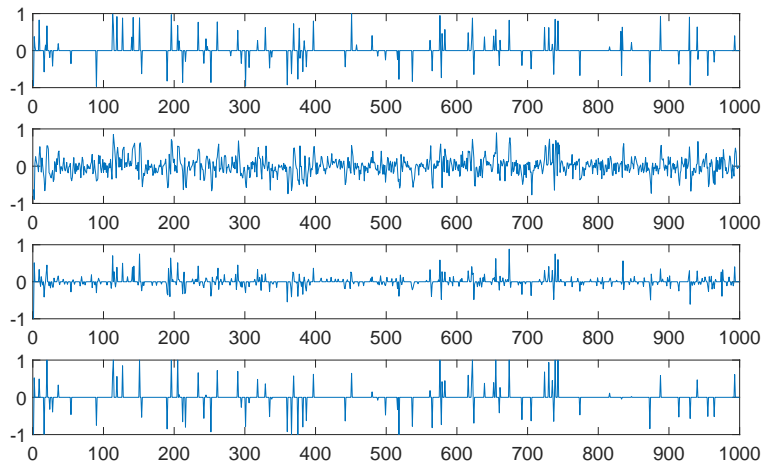}}
  \caption{Reconstruction of higher-dimensional signals ($T=1000$).
    From top to bottom: the original signal $\bxini$, the observed signal $\*y$,
    and respectively the signal reconstructed with iLASSO
    $\*x_{\mathrm{iLASSO}}$ and with our method $\hat{\*x}_{3}$.}
  ~\label{fig:ex_overlap_recons}
\end{figure}

Although our method provides good reconstruction results for medium-size
signals, handling higher-dimensional signals is highly demanding in terms of
computations as shown in our study of Section~\ref{ssec:ben_subsampl} and in
Table~\ref{table:comput_time}.
Moreover, we observed that the memory requirements of the SDP solver for its
internal process become too important.
To tackle these issues, we split the signal into smaller overlapping chunks that
are processed independently and then reassembled together.
We illustrate the example of Figure~\ref{fig:ex_overlap_recons} where we
reconstruct a signal of dimension $T=1000$ using $11$ chunks of length $100$
with $10$ overlapping samples on both extremities.
The overlapping sections are averaged in order to obtain the final signal.
The decimation operator is $D_{\infty}$ and the relaxation order is set to $3$.
We observe that our method yields a better reconstruction than iLASSO with a
mean square error of $0.43$ against $0.69$ for iLASSO\@.


\section{Conclusion}
\label{sec:conclusion}

We have proposed a method to globally solve nonconvex problems involving exact
relaxation of $\ell_{0}$ in order to reconstruct sparse signal from degraded
observations.
One of the main advantages of our method is that it is able to deal with
nonlinear degradations.
We have first reformulated our piecewise rational criterion into a rational
optimization problem before solving this problem using a hierarchy of convex SDP
relaxations that benefits from the sparsity of the rational functions.
We have then discussed the complexity of the obtained SDP and methods to
decrease both the converging relaxation order in the hierarchy and the dimension
of the SDP problem.
Finally, our simulations illustrate the domain of applicability of the method
and its high potential for finding a good approximation to a global minimum.
Although providing good results for medium-size problems, our method shows
computational limitations for larger-scale signals and filters with longer
impulse response.

\clearpage



\clearpage
\appendices

\section{Detailed computation complexity of relaxed SDP problems}
\label{sec:detail_complex_comp}

We detail here the computation of the complexity of an SDP problem, i.e.\ the
quadruple $(n,m,m_{\mathrm{s}},\ell)$, depending on the initial data like $U$, $L$
and $T$ as well as the relaxation order $k$.

\subsection{Number of blocks \texorpdfstring{$m_{\mathrm{s}}$}{ms}}
\label{ssec:nb_blocks}

In order to solve~\eqref{eq:recast_pb}, we introduce $U+T$ measures.
Moment and localizing matrices of each measure yield a block in the SDP problems
of the hierarchy.
There is one moment matrix per measure, i.e.\ a total of $U+T$ moment matrices.
The number of localizing matrices for each measure is equal to the number of
polynomial constraints defining the set $\mK_{i}$.
Equation~\eqref{eq:supp_cons_1} gives $T_{i}$ constraints for the definition of each
set $\mK_{i}$ associated to the measures ${(\mu_u)}_{u \in \nint{1,U}}$
while~\eqref{eq:supp_cons_2} gives $1+3I$ constraints for each set $\mK_{i}$
associated to the measures ${(\nu_{t})}_{t \in \nint{1,T}}$.
Indeed the polynomial equality constraint in~\eqref{eq:supp_cons_2} is translated
into two polynomial inequality constraints.
The first $L-1$ measures ${(\mu_u)}_{u \in \nint{1,U}}$ are defined on a number $T_{i}$
of variables smaller than $L$ due to the convolution filter.
In the following, we neglect it for the sake of clarity and assume that $T_{i}$ is
equal to $L$ for all the measure ${(\mu_u)}_{u \in \nint{1,U}}$.
Thus, the final number of blocks in the matrices $\*C$ and
${(\*A_i)}_{i \in \nint{1,m}}$ in~\eqref{eq:sdp_pb} is
\[
  m_{\mathrm{s}} = U(1+L) + T(2+3I) \, .
\]
It is interesting to notice that the relaxation order $k$ does not have any
effect on the number of blocks; it only increases the size of the blocks.


\subsection{Number of linear equality constraints
  \texorpdfstring{$\ell$}{l}}

We then count the number of linear equality constraints
in~\eqref{eq:sdp_eq_sum_ration}, without considering the redundant ones.
For $u$ belonging to $\nint{1,U-1}$, $\theta_{u}$ denotes the overlap parameter defined
as the number of variables shared between $g_{u}$ and $g_{u+1}$.
Note that $\theta_{u}$ depends on $u$ but also on the length of the filter $L$ and on
the parameter of the decimation $\delta$.
Furthermore, we remark that all the rational functions ${(g_{u})}_{u \in \nint{1,U}}$
have same degree at their numerator and denominator.
We denote their denominator by $q_{u}$, and we define
$\mathrm{d}_{q} = \mathrm{d}_{q_{u}}$.

Following Section~\ref{sssec:coupling}, we need to consider equality of moments
of monomials in $\theta_{u}$ variables up to degree $2(k-d_{q})$, which gives
$\binom{\theta_{u} + 2(k-d_{q})}{2(k-d_{q})}$ equality constraints for every $u$ in
$\nint{1,U-1}$ on consecutive measures ${(\mu_u)}_{u \in \nint{1,U}}$.
Adding the linear constraints linking moment related to ${(\mu_u)}_{u \in \nint{1,U}}$
and ${(\nu_{t})}_{t \in \nint{1,T}}$, we finally obtain
\[
   \ell = \sum_{u=1}^{U-1}\binom{\theta_{u} + 2(k-d_{q})}{2 (k-d_{q})} + 2(k-d_{\zeta})T
  \, ,
\]
where $d_{\zeta}$ corresponds to~\eqref{eq:defdemideg} for the maximal degree of the
denominator of rational function ${(\zeta_{i})}_{i \in \nint{1,I}}$.
The impact of linear equality constraints on the computational time of SDP
solver is minor compared to $n$, $m$ and $m_{s}$.


\subsection{Dimension of the global moment vector \texorpdfstring{$m$}{m}}

The dimension $m$ of the vector $\*v$ is simply obtained by summing up the
dimension of the moment vectors for all the measures ${(\mu_u)}_{u \in \nint{1,U}}$
and ${(\nu_{t})}_{t \in \nint{1,T}}$.
Considering ${(\mu_u)}_{u \in \nint{1,U}}$ as $U$ measures on $L$ variables and
${(\nu_{t})}_{t \in \nint{1,T}}$ as $T$ measures on $1+I$ variables, it follows that
\[
  m =  U\binom{L+2k}{2k} + T\binom{1+I+2k}{2k} \, .
\]


\subsection{Dimension of the semi-definite constraint
  \texorpdfstring{$n$}{n}}

At last, $n$ is the sum of all the block sizes of the matrices in the SDP
problem, that is the sum of the size of the all moment and localizing matrices.
The $U$ moment matrices corresponding to measures ${(\mu_u)}_{u \in \nint{1,U}}$ have
size $\binom{L+k}{k}$ while the ones corresponding to
${(\nu_{t})}_{t \in \nint{1,T}}$ have size $\binom{1+I+k}{k}$.
Since all the polynomial constraints defining the sets
${(\mK_{i})}_{i \in \nint{1,U+T}}$ are linear or quadratic, the localizing matrices
have respectively a size of $\binom{L+k-1}{k-1}$ for measures
${(\mu_u)}_{u \in \nint{1,U}}$ and $\binom{I+k}{k-1}$ for measures
${(\nu_{t})}_{t \in \nint{1,T}}$.
Finally, we obtain
\[
  n = U\left(\binom{L+k}{k} + L \binom{L+k-1}{k-1}\right)
      + T\left(\binom{1+I+k}{k} + (1+3I)\binom{I+k}{k-1} \right) \, .
\]


\end{document}